\documentclass[11pt,reqno]{amsart}
\usepackage[T1]{fontenc}
\usepackage[utf8]{inputenc}
\usepackage[left=33mm,right=33mm,top=30mm,bottom=32mm]{geometry}
\usepackage{mathtools,amssymb,amsthm,xcolor,microtype}
\usepackage[shortlabels]{enumitem}
\usepackage[normalem]{ulem}
\usepackage{needspace}
\usepackage[colorlinks,breaklinks,linkcolor=red!80!black,anchorcolor=green,citecolor=blue!60!black]{hyperref}

\definecolor{bleu1}{RGB}{0,57,128}
\def\bleu1{\color{bleu1}}

\usepackage{etoolbox}
\patchcmd{\section}{\normalfont}{\normalfont \bleu1}{}{}
\patchcmd{\subsection}{\normalfont}{\normalfont \bleu1}{}{}
\patchcmd{\subsubsection}{\normalfont}{\normalfont \bleu1}{}{}

\newtheorem{The}{\bleu1 Theorem}[section]
\newtheorem{Pro}[The]{\bleu1 Proposition}
\newtheorem{Lem}[The]{\bleu1 Lemma}
\newtheorem{Cor}[The]{\bleu1 Corollary}

\theoremstyle{definition}
\newtheorem{defn}[The]{\bleu1 Definition}
\newtheorem{Rem}[The]{\bleu1 Remark}

\def\le{\leqslant}

\def\ge{\geqslant}

\numberwithin{equation}{section}
\setlist[enumerate]{label=\textup{(\roman*)},leftmargin=*,itemsep=2pt}

\newcommand{\R}{\mathbb R}

\newcommand{\SCL}{\mathrm{SCL}}
\newcommand{\Sing}{\operatorname{Sing}}
\newcommand{\Cut}{\operatorname{Cut}}
\newcommand{\co}{\operatorname{co}}

\newcommand{\argmin}{\operatorname*{arg\,min}}
\newcommand{\argmax}{\operatorname*{arg\,max}}
\newcommand{\graph}{\operatorname{graph}}

\makeatletter
\@namedef{subjclassname@2020}{\textup{2020} Mathematics Subject Classification}
\makeatother

\hypersetup{pdftitle={Generalized Hamiltonian gradient flow of contact type II},pdfauthor={Qingyang Cai, Wei Cheng, Tianqi Shi}}

\title[Contact gradient flow: singularities and controllability]{Generalized Hamiltonian gradient flow of contact type II:\\ singularities, cut locus, and controllability}
\author{Qingyang Cai}
\address{School of Mathematics, Nanjing University, Nanjing 210093, China}
\email{652023210002@smail.nju.edu.cn}
\author{Wei Cheng}
\address{School of Mathematics, Nanjing University, Nanjing 210093, China}
\email{chengwei@nju.edu.cn}
\author{Tianqi Shi}
\address{School of Mathematics and Statistics, Northeast Normal University, Changchun 130024, China}
\email{tqshi@nenu.edu.cn}
\date{}
\begin{document}
\begin{abstract}
We study singularities, cut loci and controllability for time-dependent contact Hamilton--Jacobi equations. A contact version of Arnaud's theorem identifies the superdifferential graph of an evolved profile with the contact Hamiltonian image of the differential graph of its positive Lax--Oleinik regularization. This representation defines intrinsic singular characteristics and determines their initial velocity through a minimum-energy selection on the space-time superdifferential. The curves are continuous and locally $1/2$-H\"older, but they need not be locally Lipschitz. The exponent $1/2$ is optimal even for a quadratic Hamiltonian. The intrinsic characteristic issued from a cut point is singular at every later time in the short-time interval. We also establish the contact Lax--Oleinik commutator identities, which identify the superlevel and sublevel sets of the cut time with the corresponding terminal momentum graphs through bi-Lipschitz maps. For reachable terminal data we characterize the associated calibrated sets, and we describe all continuous initial data that reach a prescribed semiconcave profile. The proofs use Herglotz' variational principle and short-time estimates for contact fundamental solutions.
\end{abstract}
\maketitle

\section{Introduction}\label{sec:intro}

This is the second in a series of papers on generalized Hamiltonian gradient flows for contact-type Hamilton--Jacobi equations. Part~I~\cite{Cheng_Hu_Wang2026} established the regularity theory of the contact fundamental solution associated with Herglotz' variational principle. The present paper uses that theory to study singularities of viscosity solutions to time-dependent contact Hamilton--Jacobi equations, the geometry of the cut locus, and the controllability of prescribed terminal data, within the class of \emph{generalized Hamiltonian gradient flows} of \cite{CCHW2024}. Strict singular characteristics in the contact setting are treated in Part~III.

We consider the contact Hamilton--Jacobi equation
\begin{equation}\label{HJe}\tag{HJ$_e$}
 \left\{\begin{array}{l}
   \partial_tu(t,x)+H(t,x,D_xu(t,x),u(t,x))=0,\\[2pt]
   u(0,x)=\phi(x),
 \end{array}\right.\qquad x\in M,\ t>0,
\end{equation}
where $M$ is a smooth connected closed Riemannian manifold of dimension $m$ and $H$ is a time-dependent contact Tonelli Hamiltonian satisfying the standard assumptions of Section~\ref{sec:prelim}. These are strict convexity in the momentum variable, superlinear growth, boundedness of $H_u$, and a time-growth condition for $H_t$. No monotonicity assumption on the contact variable is imposed beyond the uniform bound on $H_u$.

By Herglotz' variational principle, the solution of~\eqref{HJe} is given by the negative Lax--Oleinik operator
\[
 u(t,x)=T_0^t\phi(x)=\inf_{y\in M}\{\phi(y)+h(0,t,y,x,\phi(y))\},
\]
where $h$ is the negative contact fundamental solution, defined in~\eqref{eq:h} as the infimum of the Herglotz action over curves joining $y$ to $x$ with prescribed initial scalar value. The action variable is not a free parameter but evolves along each curve through the Carath\'eodory equation
\[
 \dot u_\xi(s)=L(s,\xi(s),\dot\xi(s),u_\xi(s)).
\]
Prescribing the terminal value instead leads to the positive fundamental solution $\breve h$ and to the positive Lax--Oleinik operator $\breve T_a^b$, defined in~\eqref{eq:breveh} and \eqref{eq:LOpositive}. Every minimizer is a $C^2$ curve, and along the lifted curve with dual arc $p=L_v(\cdot,\xi,\dot\xi,u_\xi)$ the Lie equations
\[
 \dot\xi=H_p(s,\xi,p,u_\xi),\qquad
 \dot p=-H_x(s,\xi,p,u_\xi)-H_u(s,\xi,p,u_\xi)\,p,\qquad
 \dot u_\xi=\langle p,H_p\rangle-H
\]
hold. We denote the flow of this contact Hamiltonian system by $\Phi_H^{s,t}$. 
The implicit variational principle was developed in the
autonomous setting in
\cite{Su_Wang_Yan2016,Wang_Wang_Yan2017,Wang_Wang_Yan2019_1}
and extended to non-autonomous contact systems in
\cite{Wang_Yan2019}. It provides an alternative variational
framework for the representation of solutions to~\eqref{HJe}.
For the problems below we use Herglotz' principle, since the calibrated curves are defined through it. The variational construction also yields the notion of a $u$-calibrated curve (Definition~\ref{def:calibrated}), which underlies the definition of the cut locus of $u$ (Definition~\ref{def:cut}) and the analysis of the present paper. We also use the space-time reachable gradients of $u$, defined in~\eqref{eq:Ru} in terms of the reachable gradient set $D^*u(t,x)$.

The analysis of singularities of viscosity solutions of Hamilton--Jacobi equations has a long history, from the classical works on the cut locus and the singular set \cite{Albano_Cannarsa2002,Cannarsa_Sinestrari_book,Cannarsa_Yu2009,Khanin_Sobolevski2016} to the intrinsic approach of \cite{Cannarsa_Cheng3,Cannarsa_Cheng_Fathi2017,Cannarsa_Cheng_Fathi2021} for autonomous Tonelli Hamiltonians and its extension in \cite{Cannarsa_Cheng_Hong2025}. In the classical setting, Arnaud's theorem \cite{Arnaud2011} describes the evolution, under the Hamiltonian flow, of the Lagrangian submanifold associated with the gradient graph of the Lasry--Lions regularization, and relates it to the singular Lagrangian given by the superdifferential graph of the solution. The underlying $C^{1,1}$ regularization was developed in \cite{Bernard2007}. A direct proof for time-dependent non-contact Hamiltonians was given in \cite{Shi_Cheng_Hong2023}, where it served as the foundation for the study of the topology of singular sets of semiconcave functions, and in \cite{Cannarsa_Cheng_Hong2025} the Lax--Oleinik commutators were related to the topology of the singular set and to controllability. Strict singular characteristics and related variational constructions are treated in \cite{Cheng_Hong2022a,CCHW2024}, and contact $C^{1,1}$ subsolutions are studied in \cite{Wang_Xu_Zhao2026}. For the basic relations between characteristics and viscosity solutions of the contact Hamilton--Jacobi equation, including the notion of $u$-calibrated curves within Herglotz' principle, we refer to \cite{CCWY2019,CCJWY2020}. A survey of the intrinsic point of view is given in \cite{Cannarsa_Cheng2021a}.

The contact Hamiltonian $H(s,x,p,u)$ depends on the unknown function $u$ itself, which couples the state and the cost variable through the Carath\'eodory equation. This coupling produces exponential weights $e^{\int L_u\,d\tau}$ in the first-variation formulas and in the related estimates. Since $L_u=-H_u$ at dual points and no sign condition is imposed on $H_u$, the contact form may contract or expand conformally along the flow. Unlike autonomous Tonelli systems, the Hamiltonian is generally not conserved when it depends on time, and the exponential weights have to be kept in the variational estimates. The contact flow remains locally invertible, $(\Phi_H^{s,t})^{-1}=\Phi_H^{t,s}$, but the value evolution is not reversible in general. A calibrated curve may cease to admit a forward calibrated continuation, and this failure is the cut phenomenon studied below. The loss of regularity is therefore tracked through the quantitative estimates for the fundamental solutions from Part~I. The Herglotz problem and the two-parameter Lax--Oleinik evolution keep their basic form, but the classical arguments need modification.

\medskip
\noindent\textbf{Main results.}
\begin{itemize}
  \item \emph{Contact Arnaud's theorem (Theorem~\ref{thm:Arnaud}).} For $0<s<t$ with $t-s$ sufficiently small, $\breve T_s^tT_0^t\phi\in C^{1,1}(M)$ and the contact flow maps the graph of $D_x\breve T_s^tT_0^t\phi$ onto the superdifferential graph of $T_0^t\phi$, that is,
  \[
   \Phi_H^{s,t}\bigl(\graph(D_x\breve T_s^tT_0^t\phi)\bigr)
   =\graph(D_y^+T_0^t\phi),
  \]
  with all graphs including the scalar component. \cite{Wang_Xu_Zhao2026} established the existence of $C^{1,1}$ subsolutions in the implicit variational principle under stronger regularity assumptions.

  \item \emph{Intrinsic singular characteristics (Theorem~\ref{thm:ISC}).} The unique maximizing endpoint of the positive problem defines a curve $y^x$ through~\eqref{eq:ISC}. The curve is continuous on $[s,t]$ and locally $1/2$-H\"older continuous on $(s,t]$, and its right derivative at $s$ is the Hamiltonian velocity $H_p(s,x,p^\#_{u,H}(s,x),u(s,x))$, where $(q^\#_{u,H},p^\#_{u,H})$ is the unique minimizer of $q+H(s,x,p,u(s,x))$ over the full space-time superdifferential $D^+u(s,x)$ (equation~\eqref{eq:minselection}). This extends the intrinsic construction of \cite{Cannarsa_Cheng3} to the time-dependent contact setting.   
   For the evolving terminal profile $u(t,\cdot)$ the curve can fail to be locally Lipschitz, and the exponent $1/2$ is optimal even for a quadratic Hamiltonian (Remark~\ref{rem:nonLip}), while for an autonomous Hamiltonian and a fixed semiconcave terminal profile the corresponding short-time curve is Lipschitz (Remark~\ref{rem:frozen}).

  \item \emph{Propagation of singularities and the cut locus (Theorem~\ref{thm:propagation}, Corollary~\ref{cor:Cutclosure}).} If $(s,x)$ is a cut point of $u$, then the intrinsic characteristic issued from $(s,x)$ is singular at every later time in the short-time interval. Consequently
  \[
   \Sing(u)\subset\Cut(u)\subset\overline{\Sing(u)}.
  \]

  \item \emph{Lax--Oleinik commutator identities (Propositions~\ref{prop:commutators}--\ref{prop:reachability}).} The scalar contact equation preserves strict order, and the endpoint inversion between $u+h$ and $u-\breve h$ gives the operator inequalities $T_a^b\breve T_a^b\phi\ge\phi$ and $\breve T_a^bT_a^b\phi\le\phi$, together with the triple composition identities $T_a^b\breve T_a^bT_a^b\phi=T_a^b\phi$ and $\breve T_a^bT_a^b\breve T_a^b\phi=\breve T_a^b\phi$. 
The equality $T_a^b\breve T_a^b\phi=\phi$ characterizes the reachability of a prescribed terminal profile through the evolution \eqref{HJe} (Proposition~\ref{prop:reachability}). If, in addition, $\phi$ is semiconcave and $\breve T_a^b\phi\in C^{1,1}(M)$, then the graph identity and regularity persist for every intermediate initial time (Proposition~\ref{prop:continuation}).

\item \emph{Cut time function and bi-Lipschitz geometry of the cut locus (Theorem~\ref{thm:topology}).} For a fixed terminal time $b$, the super-level sets $\{x\in M:\tau_u(t,x)\ge b-t\}$ of the cut time function \eqref{eq:cuttime} are bi-Lipschitz homeomorphic to the graph $G_b^*(u)$ of the space-time reachable gradients at time $b$, and, within the $C^{1,1}$ regularization interval, the complementary sub-level sets $\{x\in M:\tau_u(t,x)<b-t\}$ are bi-Lipschitz homeomorphic to the remaining supergradient graph $G_b^\#(u)$ (Propositions~\ref{prop:contact-curves} and~\ref{prop:contact-cuttime}). The maps, given by the backward projections along the contact flow, preserve Hausdorff dimension and countable $k$-rectifiability for every $k$, which relates the cut-time problem to the geometry of the terminal momentum sets (Remark~\ref{rem:dimension}).

  \item \emph{Controllability of prescribed terminal data (Corollary~\ref{cor:controls}).} For a reachable semiconcave terminal profile $\phi$, the minimal initial datum is $\breve T_a^b\phi$, every continuous control $\psi$ with $T_a^b\psi=\phi$ dominates it and agrees with it on the initial projection of the backward trajectories issued from the spatial reachable gradients of $\phi$, and these conditions are also sufficient (Propositions~\ref{prop:min-control} and~\ref{prop:all-control-sets}).
\end{itemize}

Part~III will construct the family of \emph{strict singular characteristics} (SSC) for~\eqref{HJe}. It will establish their right-differentiability, their stability under approximation of the data, and their global propagation. It will also relate these curves to the contact dynamics through a variational principle of Herglotz type. 
The commutator defects $T_a^b\breve T_a^bf-f$ and $g-\breve T_a^bT_a^bg$, together with the sublevel structure of the cut time function, are natural objects for the study of the irreversibility of the value evolution. Their possible relations with dissipativity and with thermodynamical models of contact systems are left for future work \cite{Bravetti2017,Bravetti_Cruz_Tapias2017}. Our point of view on the underlying dynamics is close to the geometric-mechanics description of contact Hamiltonian systems developed by de Le\'on and his collaborators~\cite{de_Leon_Lainz_JMP2019,De_Leon_Lainz2019,De_Leon_Lainz_Munoz_Lecanda2023}. We return to these problems in Part~III.

The paper is organized as follows. Section~\ref{sec:prelim} collects the preliminaries on time-dependent contact Hamiltonians, Herglotz' variational principle, semiconcavity, and the theory of $u$-calibrated curves for~\eqref{HJe}. Section~\ref{sec:arnaud} contains the contact version of Arnaud's theorem, the construction and properties of the intrinsic singular characteristics, the propagation of singularities from cut points, and the Lax--Oleinik commutator identities. Section~\ref{sec:cut} is devoted to the topological structure of the cut locus and to the controllability problem. Appendix~\ref{app:A} provides the short-time estimates for the contact fundamental solutions and the smooth terminal constructions, and Appendix~\ref{app:B} collects the proofs of the calibration results used in the preliminaries.

\medskip
\noindent\textbf{Acknowledgments.}
Wei Cheng was supported in part by the National Natural Science Foundation of China (Grant No.~12231010). Tianqi Shi was supported in part by the National Natural Science Foundation of China (Grant No.~12601348).

\smallskip
\noindent\textbf{Declaration of generative AI use.}
During the preparation of this work, the authors used DeepSeek and ChatGPT for language editing and formatting. After using these tools, the authors reviewed and edited the content and take full responsibility for the content of the publication.

\section{Preliminaries}\label{sec:prelim}
\subsection{Herglotz's variational principle and contact fundamental solutions}
We work with $L\in C^2(\R\times TM\times\R)$ and its Legendre dual
\begin{equation}\label{eq:duality}
 H(t,x,p,u)=\sup_{v\in T_xM}\{\langle p,v\rangle-L(t,x,v,u)\}.
\end{equation}
We assume that $H$ is also $C^2$. On each compact time interval the hypotheses below hold, with constants uniform in time.
\begin{enumerate}[label=\textup{(L\arabic*)}]
\item $L(t,x,\cdot,u)$ is strictly convex.
\item There are superlinear functions $\theta_0,\theta_1$ and constants $c_0,c_1$ such that
$\theta_0(|v|)-c_0\le L(t,x,v,0)\le\theta_1(|v|)+c_1$.
\item $|L_u(t,x,v,u)|\le K$.
\item $|L_t(t,x,v,u)|\le C_1+C_2L(t,x,v,u)$ for some $C_1,C_2\ge0$.
\end{enumerate}
The dual hypotheses are the following, with the corresponding constants.
\begin{enumerate}[label=\textup{(H\arabic*)}]
\item $H(t,x,\cdot,u)$ is strictly convex.
\item $\theta_1^*(|p|)-c_1\le H(t,x,p,0)\le\theta_0^*(|p|)+c_0$, where $\theta_i^*$ denotes the Fenchel--Legendre dual of $\theta_i$.
\item $|H_u(t,x,p,u)|\le K^*$ for some $K^*>0$.
\item There are $C_1^*,C_2^*\ge0$ such that
\[
 |H_t(t,x,p,u)|\le C_1^*+C_2^*\{\langle p,H_p(t,x,p,u)\rangle-H(t,x,p,u)\}.
\]
\end{enumerate}
At dual points, $H_u=-L_u$ and $H_t=-L_t$, so the constants can be chosen consistently. We use $K$ for a common scalar Lipschitz bound when applying comparison.
The $C^2$ Legendre duality implies
$L_{vv}H_{pp}=I$ at dual points, so both Hessians are positive definite. This nondegeneracy is used in the short-time estimates. We write $d$ for the Riemannian distance and $\pi_x:T^*M\times\R\to M$ for the spatial projection. Constants below may depend on a fixed compact time interval.

For $s<t$ and $x,y\in M$, retain the endpoint class
\[
 \Gamma_{x,y}^{s,t}:=\{\xi\in W^{1,1}([s,t],M):\xi(s)=x,\ \xi(t)=y\}.
\]
For $\xi\in\Gamma_{x,y}^{s,t}$, consider the initial and terminal Carath\'eodory problems
\begin{align}
 \dot u_\xi(\tau)&=L(\tau,\xi(\tau),\dot\xi(\tau),u_\xi(\tau)),
                       &u_\xi(s)&=u,\label{eq:Caratheodory}\\
 \dot w_\xi(\tau)&=L(\tau,\xi(\tau),\dot\xi(\tau),w_\xi(\tau)),
                       &w_\xi(t)&=u.\label{eq:terminalCaratheodory}
\end{align}
The admissible classes are
\begin{align*}
 \Gamma_{x,y}^{s,t}(u)&:=\{\xi\in\Gamma_{x,y}^{s,t}:
                  \eqref{eq:Caratheodory}\text{ has an absolutely continuous solution}\},\\
 \breve\Gamma_{x,y}^{s,t}(u)&:=\{\xi\in\Gamma_{x,y}^{s,t}:
                  \eqref{eq:terminalCaratheodory}\text{ has an absolutely continuous solution}\}.
\end{align*}
By (L3), both classes consist precisely of curves with
$L(\cdot,\xi,\dot\xi,0)\in L^1(s,t)$, so they are independent of the scalar endpoint value. Indeed, an absolutely continuous scalar solution is bounded, and (L3) makes integrability at that solution equivalent to integrability at zero. Conversely, integrability at zero and the scalar Lipschitz bound give global existence and uniqueness on the prescribed finite interval.

The negative and positive contact fundamental solutions are
\begin{align}
 h(s,t,x,y,u)&:=\inf_{\xi\in\Gamma_{x,y}^{s,t}(u)}
       \int_s^tL(\tau,\xi(\tau),\dot\xi(\tau),u_\xi(\tau))\,d\tau,\label{eq:h}\\
 \breve h(s,t,x,y,u)&:=\inf_{\xi\in\breve\Gamma_{x,y}^{s,t}(u)}
       \int_s^tL(\tau,\xi(\tau),\dot\xi(\tau),w_\xi(\tau))\,d\tau.\label{eq:breveh}
\end{align}
Thus $u+h$ is the least terminal value for initial value $u$, whereas $u-\breve h$ is the greatest initial value for terminal value $u$. Along a specified lift we also write $u_\xi$ for the scalar component when its endpoint convention has already been fixed.

\begin{Pro}[Herglotz's variational principle]\label{prop:Herglotz}
Under the preceding hypotheses, the infima in \eqref{eq:h} and \eqref{eq:breveh} are attained. Every minimizing curve $\xi$ and its scalar lift $u_\xi$ are $C^2$. With
$p_\xi=L_v(\tau,\xi,\dot\xi,u_\xi)$, the lifted curve solves
\begin{equation*}
\tag{Lie}\label{eq:contact}
 \dot\xi=H_p,\qquad \dot p_\xi=-H_x-H_up_\xi,\qquad
 \dot u_\xi=\langle p_\xi,H_p\rangle-H.
\end{equation*}
Equivalently, the dual arc satisfies the Herglotz Euler--Lagrange equation
\[
 \frac d{d\tau}L_v(\tau,\xi,\dot\xi,u_\xi)
 =L_x(\tau,\xi,\dot\xi,u_\xi)
           +L_u(\tau,\xi,\dot\xi,u_\xi)L_v(\tau,\xi,\dot\xi,u_\xi).
\]
These statements hold for both endpoint conventions.
\end{Pro}
The existence, regularity, and a priori estimates are the Herglotz variational results of \cite{CCWY2019,CCJWY2020}. The positive problem also follows by the reversal in Lemma~\ref{lem:reversal}. We denote the flow of \eqref{eq:contact} by $\Phi_H^{s,t}$. We need this flow only along the relevant compact families, and we do not assume that it is complete on the whole phase space.

The negative and positive Lax--Oleinik operators are
\begin{align}
 T_s^t\phi(x)&=\inf_{y\in M}\{\phi(y)+h(s,t,y,x,\phi(y))\},\label{eq:LOnegative}\\
 \breve T_s^t\phi(x)&=\sup_{y\in M}\{\phi(y)-\breve h(s,t,x,y,\phi(y))\}.\label{eq:LOpositive}
\end{align}
Both preserve order, and for $a\le s\le b$ they satisfy
\begin{equation}\label{eq:Markov}
 T_a^b=T_s^bT_a^s,\qquad
 \breve T_a^b=\breve T_a^s\breve T_s^b.
\end{equation}
The operators at equal times are the identity. Their sup-norm stability estimate is
\begin{equation}\label{eq:stability}
 \|T_a^bf-T_a^bg\|_\infty,
 \ \|\breve T_a^bf-\breve T_a^bg\|_\infty
 \le e^{K(b-a)}\|f-g\|_\infty.
\end{equation}
These statements follow by concatenating the scalar equations and applying scalar comparison. A function is $(\kappa_1,\kappa_2)$-Lipschitz in the large if
$|\phi(x)-\phi(y)|\le\kappa_1d(x,y)+\kappa_2$. For lower semicontinuous such data, $u(t,x)=T_0^t\phi(x)$ is locally Lipschitz and locally semiconcave on $(0,\infty)\times M$, with uniform bounds on compact positive time strips, and solves \eqref{HJe}. When the initial datum is continuous, the solution is continuous up to time zero. Positive evolutions of continuous terminal data are locally semiconvex jointly in time and space away from the terminal time. To see the joint statement, fix $b$ and a continuous terminal datum $f$, and set $\widetilde L_b(r,x,v,w)=L(b-r,x,-v,-w)$. Let $\widetilde T$ be the negative Lax--Oleinik operator associated with $\widetilde L_b$. The curve and scalar reversals give
\[
 -\breve T_{b-r}^{b}f(x)=\widetilde T_0^r(-f)(x),\qquad r>0.
\]
The right side is jointly semiconcave on positive time intervals. Hence $(t,x)\mapsto\breve T_t^bf(x)$ is jointly semiconvex for $t<b$.

\begin{Pro}\label{prop:domination}
Let $u$ be a locally semiconcave viscosity subsolution of \eqref{HJe}. For $0<t_1<t_2$ and $x,y\in M$,
\begin{align}
 u(t_2,y)&\le u(t_1,x)+h(t_1,t_2,x,y,u(t_1,x)),\label{eq:domh}\\
 u(t_1,x)&\ge u(t_2,y)-\breve h(t_1,t_2,x,y,u(t_2,y)).\label{eq:dombh}
\end{align}
For $u=T_0^t\phi$ with lower semicontinuous initial data as above, both inequalities also hold when $t_1=0$.
\end{Pro}
The proof, including the initial-time assertion, is given in Appendix~\ref{app:B}.

\subsection{Semiconcavity and calibrated curves}
Write $\SCL(M)$ for the functions semiconcave with linear modulus, and $\SCL_{\mathrm{loc}}$ for its local version. For $f\in\SCL(M)$, the superdifferential and reachable-gradient set are
\begin{align*}
D^+f(x)&:=\{D\psi(x):\psi\in C^1\text{ near }x,\ \psi\ge f,\ \psi(x)=f(x)\},\\
D^*f(x)&:=\{\lim_{j\to\infty}Df(x_j):x_j\to x,\ x_j\notin\Sing(f)\}.
\end{align*}
We write $\Sing(u)$ for the set where $u$ is not differentiable in space-time.
For a space-time function $u$, the notation $D_x^+u(t,x)$ and $D_x^*u(t,x)$ refers to the function $u(t,\cdot)$. We recall the standard identities
\begin{equation}\label{eq:SCfacts}
D^+f=\co D^*f,\qquad
D_x^+u(t,x)=\Pi_xD^+u(t,x),\qquad
d^+f(z;v)=\min_{p\in D^+f(z)}\langle p,v\rangle,
\end{equation}
where $\Pi_x(q,p)=p$ (see, e.g., \cite[Chapter~3]{Cannarsa_Sinestrari_book}, and in particular \cite[Lemma~3.3.16]{Cannarsa_Sinestrari_book} for the projection formula). The second equality must not be transferred to reachable gradients without further justification. We use the convention that the quadratic remainder in the semiconcavity inequality is $C|h|^2/2$.

\begin{Pro}\label{prop:marginal}
Let $S$ be compact and $F:S\times\Omega\to\R$ be continuous, where $\Omega\subset\R^n$ is open. Suppose $F(\sigma,\cdot)$ is $C^2$ with locally uniform $C^2$ bounds. Then $f(x)=\min_{\sigma\in S}F(\sigma,x)$ is locally semiconcave. Writing $M(x)=\argmin_{\sigma\in S}F(\sigma,x)$, one has
\[
D^*f(x)\subset\{D_xF(\sigma,x):\sigma\in M(x)\},\qquad
D^+f(x)=\co\{D_xF(\sigma,x):\sigma\in M(x)\}.
\]
In particular, a unique active gradient implies differentiability.
\end{Pro}
This is the standard marginal-function characterization of semiconcavity (see \cite[Section~3.4]{Cannarsa_Sinestrari_book}).

\begin{Lem}\label{lem:SCclosed}
Suppose $f_j\to f$ locally uniformly with locally uniform semiconcavity constants. If $x_j\to x$ and $p_j\in D^+f_j(x_j)$, then $\operatorname{dist}(p_j,D^+f(x))\to0$.
\end{Lem}
This is the local version of \cite[Lemma~2.2]{CCHW2024}. The manifold version follows in local charts.

\begin{Pro}\label{prop:touch}
Suppose $u$ is $C$-semiconcave, $v$ is $C$-semiconvex on $\R^n$, and $u\ge v$. At every point of $E=\{u=v\}$ both functions are differentiable and $Du=Dv$. This common gradient is $4C$-Lipschitz on $E$. The corresponding local assertion holds in charts, including space-time charts.
\end{Pro}
This is \cite[Proposition~2.5]{Cannarsa_Cheng_Hong2025}. The local version follows by restricting to a smaller coordinate neighborhood.

\begin{defn}\label{def:calibrated}
A $C^1$ curve $\eta:I\to M$ is $u$-calibrated if, for every $[a,b]\subset I$,
\[
 u(b,\eta(b))-u(a,\eta(a))
 =\int_a^bL(s,\eta(s),\dot\eta(s),u_\eta(s))\,ds,
 \qquad u_\eta(a)=u(a,\eta(a)),
\]
where $u_\eta$ solves \eqref{eq:Caratheodory} along $\eta$.
\end{defn}

\begin{Pro}\label{prop:cal-equivalence}
Let $0<a<b$, let $u$ be a locally semiconcave viscosity solution, and let $\eta:[a,b]\to M$ be admissible and absolutely continuous. The following conditions are equivalent.
\begin{enumerate}
\item $\eta$ is $u$-calibrated.
\item $u(b,\eta(b))=u_\eta(b)$, where $u_\eta(a)=u(a,\eta(a))$.
\item $u(b,\eta(b))-u(a,\eta(a))=\int_a^bL(s,\eta,\dot\eta,u(s,\eta(s)))\,ds$.
\end{enumerate}
In this case $u_\eta(s)=u(s,\eta(s))$ throughout $[a,b]$, and $\eta$ is $C^2$.
\end{Pro}

\begin{Lem}\label{lem:cal-existence}
Let $u(t,x)=T_0^t\phi(x)$ for lower semicontinuous data Lipschitz in the large. Every $(t,x)$ with $t>0$ has an incoming calibrated curve on $[0,t]$. Moreover, for $0\le t<\tau$,
\begin{equation}\label{eq:cal-contact}
 u(t,x)=\breve T_t^\tau[u(\tau,\cdot)](x)
\end{equation}
if and only if a calibrated curve on $[t,\tau]$ starts at $x$.
\end{Lem}
\begin{proof}
An attained minimizer for $T_0^t\phi(x)$ is calibrated by dynamic programming. For \eqref{eq:cal-contact}, a calibrated curve gives equality in the domination inequality. Conversely, let $\eta$ be a positive optimizer. Equality in \eqref{eq:cal-contact} implies that its terminal-value scalar lift satisfies $w_\eta(t)=u(t,x)$. By uniqueness of the scalar equation, $w_\eta$ coincides with the scalar lift issued from $u(t,x)$, and hence Proposition~\ref{prop:cal-equivalence} shows that $\eta$ is calibrated. At $t=0$ we argue directly. Let $\eta$ be the positive optimizer and $w$ its scalar lift, so that $w(0)=u(0,x)=\phi(x)$ and $w(\tau)=u(\tau,\eta(\tau))$. For $0<r<\tau$, the initial-time part of Proposition~\ref{prop:domination} gives $u(r,\eta(r))\le w(r)$. Let $E_{r,\tau}$ be the scalar endpoint map along $\eta|_{[r,\tau]}$. Domination at positive times and monotonicity of $E_{r,\tau}$ give
\[
 u(\tau,\eta(\tau))\le E_{r,\tau}(u(r,\eta(r)))\le E_{r,\tau}(w(r))=u(\tau,\eta(\tau)).
\]
Since $E_{r,\tau}$ is strictly increasing, $u(r,\eta(r))=w(r)$ for every $r\in(0,\tau)$. Together with the equalities at the endpoints, this is the calibration condition on $[0,\tau]$.
\end{proof}

\begin{defn}\label{def:cut}
A point $(t,x)\in(0,\infty)\times M$ is a cut point of $u$ if no calibrated curve starting there extends to a later time. The cut locus is
\[
 \Cut(u)=\{(t,x):u(t,x)\ne\breve T_t^\tau[u(\tau,\cdot)](x)
                          \text{ for every }\tau>t\}.
\]
\end{defn}

\Needspace{8\baselineskip}
\begin{The}\label{thm:cal-lifts}
Let $u=T_0^t\phi$ be as in Lemma~\ref{lem:cal-existence}.
\begin{enumerate}
\item If $\gamma:[0,t]\to M$ is calibrated and $p=L_v(s,\gamma,\dot\gamma,u(s,\gamma))$, then $u$ is differentiable at $(s,\gamma(s))$ for $0<s<t$, and
\[
 Du(s,\gamma(s))=(-H(s,\gamma(s),p(s),u(s,\gamma(s))),p(s)).
\]
At the terminal point $p(t)\in D_x^+u(t,\gamma(t))$.
\item There is a bijection between $D^*u(t,x)$ and the calibrated lifts ending at $(t,x)$, given by
\begin{equation}\label{eq:reachableflow}
 (\gamma(s),p(s),u_\gamma(s))=\Phi_H^{t,s}(x,p,u(t,x)),
 \qquad (q,p)\in D^*u(t,x),
\end{equation}
where $q=-H(t,x,p,u(t,x))$.
\end{enumerate}
\end{The}
The proofs of Proposition~\ref{prop:cal-equivalence} and Theorem~\ref{thm:cal-lifts} are in Appendix~\ref{app:B}. In particular, the terminal momenta of all incoming calibrated curves form
\begin{equation}\label{eq:Ru}
 \mathcal R_u(t,x):=\Pi_xD^*u(t,x)
 =\{p:(-H(t,x,p,u(t,x)),p)\in D^*u(t,x)\}.
\end{equation}
One has $D_x^*u(t,x)\subset\mathcal R_u(t,x)\subset D_x^+u(t,x)$. Indeed, spatial differentiability implies full differentiability. All reachable space-time gradients have the same spatial component and, by the equation, the same time component. Hence the first inclusion follows by approximating an element of $D_x^*u(t,x)$ by spatial gradients at differentiability points and passing to the limit. The second inclusion follows from $D^*u(t,x)\subset D^+u(t,x)$ and the projection identity in \eqref{eq:SCfacts}.

\begin{Cor}\label{cor:SingCut}
For $u=T_0^t\phi$, one has $\Sing(u)\subset\Cut(u)$.
\end{Cor}
\begin{proof}
A forward calibrated extension concatenates with an incoming calibrated curve. The junction is an interior point of a calibrated curve, so Theorem~\ref{thm:cal-lifts} gives differentiability there.
\end{proof}

\section{Arnaud's theorem and intrinsic singular characteristics}\label{sec:arnaud}
Throughout this section, $\phi$ is lower semicontinuous and Lipschitz in the large, $u(t,x)=T_0^t\phi(x)$, and $0<\tau_1<\tau_2$. The short-time constants are uniform for $s,t\in[\tau_1,\tau_2]$. All graph notation includes the contact component. For example,
\[
 \graph(D^+f)=\{(x,p,f(x)):x\in M,\ p\in D^+f(x)\}.
\]
The short-time properties of $\breve h$ used below are proved in Appendix~\ref{app:A}.

\subsection{A contact version of Arnaud's theorem}

\begin{Lem}\label{lem:twist}
Fix $R,S>0$ and write $(X(s),P(s),U(s))=\Phi_H^{t,s}(y,p,u)$. There is $t_{R,S}>0$, uniform in $y\in M$ and $s,t\in[\tau_1,\tau_2]$, such that for $0<t-s\le t_{R,S}$ the map
\[
 \Phi_{s,t,y}:B_y(0,R)\times(-S,S)\longrightarrow M\times\R,
 \qquad(p,u)\longmapsto(X(s),U(s))
\]
is a $C^1$ diffeomorphism onto its image. Here $B_y(0,R)\subset T_y^*M$. For each fixed $u\in(-S,S)$, the map
\[
 \Psi_{s,t,y,u}:B_y(0,R)\longrightarrow M,\qquad p\longmapsto X(s)
\]
is also a $C^1$ diffeomorphism onto its image.
\end{Lem}
\begin{proof}
Put $h=t-s$. In a common chart, differentiating the flow gives uniformly
\begin{equation}\label{eq:twist-estimates}
 X_p(s)=-hH_{pp}(t,y,p,u)+o(h),\quad
 X_u(s)=O(h),\quad U_p(s)=O(h),\quad U_u(s)=1+O(h).
\end{equation}
Consequently the symmetric part of $X_p$ is bounded above by $-chI$. Integration along a segment in a convex bounded momentum set proves injectivity for fixed $u$ and gives
$|X(s,p_1,u)-X(s,p_2,u)|\ge ch|p_1-p_2|$.

For joint injectivity, equal values of $U(s)$ imply
$|u_1-u_2|\le Ch|p_1-p_2|$. Equal values of $X(s)$ and the preceding monotonicity estimate then imply
$ch|p_1-p_2|\le Ch|u_1-u_2|\le Ch^2|p_1-p_2|$.
Thus both differences vanish when $h$ is small. The contact variation identity gives
\[
 U_p=P(s)X_p,\qquad
 U_u-P(s)X_u=\exp\!\left(\int_s^t H_u(\tau,X,P,U)\,d\tau\right).
\]
Subtracting $P(s)$ times the first block row in the Jacobian shows that its determinant is the positive exponential factor times $\det X_p$. The inverse function theorem completes the proof. These arguments are local in the terminal base point and uniform in a finite atlas.
\end{proof}

\begin{The}[Arnaud theorem of contact type]\label{thm:Arnaud}
There are positive constants $t_\phi\le t'_\phi\le t''_\phi$ such that, whenever $\tau_1\le s<t\le\tau_2$, the following statements hold.
\begin{enumerate}
\item if $t-s\le t''_\phi$, then
\[
 \Phi_H^{s,t}\bigl(\graph(D_x^*\breve T_s^tT_0^t\phi)\bigr)
 \subset\graph(D_y^+T_0^t\phi).
\]
\item if $t-s\le t'_\phi$, then $\breve T_s^tT_0^t\phi\in C^{1,1}(M)$, so its reachable gradient is its differential.
\item if $t-s\le t_\phi$, then
\begin{equation}\label{eq:Arnaud}
 \Phi_H^{s,t}\bigl(\graph(D_x\breve T_s^tT_0^t\phi)\bigr)
 =\graph(D_y^+T_0^t\phi).
\end{equation}
\end{enumerate}
For families of terminal profiles with uniform bounds on the norm,
Lipschitz constant, and semiconcavity constant, the same conclusions
hold uniformly with the time endpoints in any fixed compact time
interval $I$. Lemma~\ref{lem:twist} also holds uniformly for $s,t\in I$.
\end{The}
\begin{proof}
Write $f=u(t,\cdot)$, $h=t-s$, and
\[
 F_x(y)=f(y)-\breve h(s,t,x,y,f(y)),\qquad
 v(x)=\breve T_s^tf(x).
\]
Choose the short-time cone large enough to contain all maximizing
endpoints and the short backward trajectories issued from
$\graph(D^+f)$. By Lemma~\ref{lem:localization} and
Proposition~\ref{prop:kernel}, choose $t_\phi''>0$, uniform for
$\tau_1\le s<t\le\tau_2$, so that the required endpoint estimates
hold on this cone whenever $0<h\le t_\phi''$. Fix $p\in D_x^*v(x)$ and choose differentiability points $x_j\to x$ with $Dv(x_j)\to p$. Pick a maximizer $y_j$ of $F_{x_j}$. Localization and compactness allow a subsequence with $y_j\to y$. Write $P_j(\tau)$ for the momentum of the contact lift joining $(x_j,v(x_j))$ to $(y_j,f(y_j))$. Near $x_j$, the function
\[
 z\longmapsto f(y_j)-\breve h(s,t,z,y_j,f(y_j))
\]
is a smooth lower support of $v$ touching at $x_j$. Hence
\[
 Dv(x_j)=-\breve h_x(s,t,x_j,y_j,f(y_j))=P_j(s).
\]
At the maximizing endpoint, Proposition~\ref{prop:chain} applied to $(z,u)\mapsto u-\breve h(s,t,x,z,u)$ and Fermat's rule give
\[
 0\in-\breve h_y(s,t,x_j,y_j,f(y_j))
 +(1-\breve h_u(s,t,x_j,y_j,f(y_j)))D^+f(y_j).
\]
Since $\breve h_y=\mu_jP_j(t)$ and $1-\breve h_u=\mu_j>0$, we obtain $P_j(t)\in D^+f(y_j)$. Continuous dependence of the contact flow and closedness of the graph of $D^+f$ (Lemma~\ref{lem:SCclosed}) pass this relation to the limit. The terminal state of $\Phi_H^{s,t}(x,p,v(x))$ is therefore $(y,P(t),f(y))$ with $P(t)=\lim_jP_j(t)\in D^+f(y)$, hence
\[
 \Phi_H^{s,t}\bigl(\graph(D_x^*v)\bigr)\subset\graph(D^+f),
\]
proving~(1).

For~(2), decrease the threshold and choose $t_\phi'\le t_\phi''$ so that
Lemma~\ref{lem:strictconcavity} and the corresponding mixed derivative
estimates hold for $0<h\le t_\phi'$. Then $F_x$ has a unique global
maximizer, denoted $Y(x)$, and compactness gives continuity of $Y$.

For nearby $x_1,x_2$, set $y_i=Y(x_i)$. The quadratic drop at a maximizer
of a strongly concave function implies
\[
 \frac c h|y_1-y_2|^2
 \le [F_{x_1}(y_1)-F_{x_2}(y_1)]
    -[F_{x_1}(y_2)-F_{x_2}(y_2)].
\]
The right side is at most
$Ch^{-1}|x_1-x_2||y_1-y_2|$. Indeed, integrate the $x$ derivative
between $x_1,x_2$ and use
$|\breve h_{xy}|\le C/h$, $|\breve h_{xu}|\le C$, and the Lipschitz
constant of $f$. Thus $Y$ is locally Lipschitz. Differentiating a maximum
with a unique optimizer gives
\begin{equation}\label{eq:vgradient}
 Dv(x)=-\breve h_x(s,t,x,Y(x),f(Y(x)))=P(s).
\end{equation}
The endpoint derivatives are locally Lipschitz on each fixed short-time
cone. Hence $Dv$ is Lipschitz and $v\in C^{1,1}$. Thus~(2) holds for
$0<h\le t_\phi'$.

Finally, decrease the threshold once more and choose
$t_\phi\le t_\phi'$ so that Lemma~\ref{lem:twist} applies for
$0<h\le t_\phi$. Conversely, take
$(y,p,f(y))\in\graph(D^+f)$ and follow its contact trajectory backward
to $x=X(s)$. Its velocity and all phase variables remain uniformly
bounded. The fixed-endpoint minimizer defining
$\breve h(s,t,x,y,f(y))$ has bounded terminal momentum as well.
Both trajectories have the same terminal base point and scalar value,
as well as the same initial base point. Hence Lemma~\ref{lem:twist}
identifies their terminal momenta and thus the trajectories.
Consequently, $0\in D^+F_x(y)$. Strong concavity makes $y$ the unique
maximizer on the cone, and localization makes it the global maximizer.
Equation~\eqref{eq:vgradient} gives the reverse inclusion in
\eqref{eq:Arnaud}, proving~(3).

All constants can be chosen uniformly on the fixed positive compact
time strip.
\end{proof}

\subsection{Intrinsic singular characteristics}
For $(s,x)\in(0,\infty)\times M$, define
\begin{equation}\label{eq:minselection}
 (q^\#_{u,H},p^\#_{u,H})(s,x)
 :=\argmin_{(q,p)\in D^+u(s,x)}
                 \{q+H(s,x,p,u(s,x))\}.
\end{equation}
The minimizer exists by compactness and is unique. Strict convexity determines $p$, after which the minimum determines $q$.

\begin{Lem}\label{lem:selection}
At every fixed $(s,x)$ with $s>0$,
\begin{equation}\label{eq:selectionlimit}
 \lim_{t\downarrow s}D_x[\breve T_s^tT_0^t\phi](x)
 =p^\#_{u,H}(s,x)\in D_x^+u(s,x).
\end{equation}
\end{Lem}
\begin{proof}
Let $(X_h,P_h,U_h)$ be the positive optimizing lift from $x$ at time $s$ to $y_h$ at time $s+h$. The preceding theorem gives
\[
 P_h(s)=D_x\breve T_s^{s+h}u(s+h,\cdot)(x),\qquad
 P_h(s+h)\in D_x^+u(s+h,y_h).
\]
Localization and the a priori estimates keep these lifts in a compact set. Moreover, $U_h(s)=u(s,x)+O(h)$ and $P_h(s+h)-P_h(s)=O(h)$.
Choose a subsequence along which both momenta converge to $\bar p$. By \eqref{eq:SCfacts}, choose $q_h$ with
$(q_h,P_h(s+h))\in D^+u(s+h,y_h)$. Passing to a further subsequence and using Lemma~\ref{lem:SCclosed} gives $(\bar q,\bar p)\in D^+u(s,x)$. The flow implies
\[
 \frac{y_h-x}{h}\longrightarrow
 \bar v:=H_p(s,x,\bar p,u(s,x)).
\]
Semiconcave functions have Hadamard directional derivatives. Therefore
\[
 \frac{u(s+h,y_h)-u(s,x)}h
 \longrightarrow\min_{(q,p)\in D^+u(s,x)}(q+\langle p,\bar v\rangle).
\]
The supergradient inequality at $(s+h,y_h)$, evaluated at $(s,x)$, says that this limit is at least $\bar q+\langle\bar p,\bar v\rangle$. The reverse inequality follows from membership of $(\bar q,\bar p)$ in $D^+u(s,x)$. Hence
\[
 q-\bar q+\langle\bar v,p-\bar p\rangle\ge0
 \qquad\text{for all }(q,p)\in D^+u(s,x).
\]
Convexity of $H$ then gives
$q+H(s,x,p,u)-\bar q-H(s,x,\bar p,u)\ge0$.
Thus $(\bar q,\bar p)$ is the unique minimizer in \eqref{eq:minselection}. Every convergent subsequence has the same limit, proving \eqref{eq:selectionlimit}.
\end{proof}

For $0<t-s\le t_\phi$ define the intrinsic curve based at $(s,x)$ by
\begin{equation}\label{eq:ISC}
 y^x(s)=x,\qquad
 y^x(\tau)=\argmax_{y\in M}
 \{T_0^\tau\phi(y)-\breve h(s,\tau,x,y,T_0^\tau\phi(y))\},
 \quad s<\tau\le t.
\end{equation}
At cut points these curves are intrinsic singular characteristics, as justified below.

\begin{The}\label{thm:ISC}
The curve \eqref{eq:ISC} has the properties below.
\begin{enumerate}
\item for $s<\tau\le t$,
\[
 y^x(\tau)=\pi_x\Phi_H^{s,\tau}
 \bigl(x,D_x\breve T_s^\tau T_0^\tau\phi(x),
                       \breve T_s^\tau T_0^\tau\phi(x)\bigr).
\]
\item $y^x\in C([s,t],M)\cap C^{0,1/2}_{\mathrm{loc}}((s,t],M)$.
\item its right derivative at $s$ exists and equals
\[
 \dot y^x(s+)=H_p(s,x,p^\#_{u,H}(s,x),u(s,x)).
\]
\end{enumerate}
\end{The}
\begin{proof}
The flow representation follows from Theorem~\ref{thm:Arnaud}. Uniqueness of the maximizer and compactness give continuity for $\tau>s$, while localization gives $d(y^x(\tau),x)\le C(\tau-s)$ and hence continuity at $s$.

For the H\"older estimate, put $F_\tau(z)=u(\tau,z)-\breve h(s,\tau,x,z,u(\tau,z))$ and fix $J=[s+\delta,t]\Subset(s,t]$. By localization, after choosing $\Lambda>\lambda$, the two maximizers $y^x(\tau)$ and $y^x(\sigma)$ belong to a common $\Lambda$-cone whenever $\tau,\sigma\in J$ are sufficiently close. Decreasing $t_\phi$ if necessary, Proposition~\ref{prop:kernel} applies uniformly on this cone. Hence the endpoint derivative formulas and the time-Lipschitz continuity of $u$ on $J\times M$ give
\[
|F_\tau(z)-F_\sigma(z)|\le C_J|\tau-\sigma|,
\qquad z\in\{y^x(\tau),y^x(\sigma)\}.
\]
Since the Lipschitz and semiconcavity constants of $u(r,\cdot)$ are uniform for $r\in J$, Lemma~\ref{lem:strictconcavity} applies uniformly. Together with maximality,
\[
\frac{c}{\tau-s}|y^x(\tau)-y^x(\sigma)|^2
\le F_\tau(y^x(\tau))-F_\tau(y^x(\sigma))
\le 2C_J|\tau-\sigma|.
\]
This proves local $1/2$-H\"older continuity.

Finally, let $(X_h,P_h,U_h)$ be the optimizing lift on $[s,s+h]$. Then
\[
\frac{y^x(s+h)-x}{h}
=\frac1h\int_s^{s+h}H_p(r,X_h(r),P_h(r),U_h(r))\,dr.
\]
By Lemma~\ref{lem:selection} and continuous dependence of the flow, the right-hand side converges to $H_p(s,x,p^\#_{u,H}(s,x),u(s,x))$ as $h\downarrow0$. This proves the asserted right derivative.
\end{proof}

\begin{Rem}\label{rem:nonLip}
Even for $H(p)=|p|^2/2$ on $\R^2$, a globally Lipschitz initial datum can produce an intrinsic curve
\begin{equation}\label{eq:counter-curve}
 y(t)=\argmax_z\left\{u(t,z)-\frac{|z|^2}{2(t-s)}\right\},
 \qquad u(t)=T_0^t\phi,
\end{equation}
that is not locally Lipschitz at an interior time. More strongly, it fails to be locally $\alpha$-H\"older for every $\alpha>1/2$, proving that the exponent $1/2$ in Theorem~\ref{thm:ISC} is optimal. The example can also be placed on a flat two-torus, within an arbitrarily short interval after $s$.

Fix $A>0$, put $d_n=2^{-n}$, $e_n=An(n+1)2^{-n}$, and define, for $n\ge1$,
\begin{equation}\label{eq:counter-sequences}
 \begin{gathered}
 a_1=b_1=0,\quad \lambda_1=A,\quad c_1^2=8,\\
 a_{n+1}-a_n=d_n,\quad \lambda_{n+1}-\lambda_n=e_n,\\
 b_{n+1}-b_n=-1/[An(n+1)],\quad
 c_{n+1}^2-c_n^2=2\lambda_n(b_{n+1}-b_n).
 \end{gathered}
\end{equation}
Then $a_n\uparrow1$, $b_n\downarrow-1/A$, $\lambda_n\uparrow\lambda_*=9A$ and $c_n^2\downarrow2$. The last identity follows from
$\sum_{n\ge1}\lambda_n/[An(n+1)]=1+\sum_{n\ge1}n2^{-n}=3$.
Put
\[
 p_{n,\pm}=(\pm c_n,a_n),\qquad
 \alpha_n=b_n+\frac s2(c_n^2+a_n^2),\qquad
 \phi(z)=\inf_{n,\pm}(p_{n,\pm}\cdot z+\alpha_n).
\]
The limiting index $n=\infty$ may be included without changing the infimum. Boundedness of the slopes makes $\phi$ globally Lipschitz and concave. Commuting two infima in the Hopf--Lax formula gives
\begin{equation}\label{eq:counter-solution}
 u(t,z)=\inf_{n,\pm}\left\{p_{n,\pm}\cdot z+\alpha_n-\frac t2|p_{n,\pm}|^2\right\}.
\end{equation}
Write $h=t-s$, $\lambda=h^{-1}$, and $z=hv$. Then
\[
 u(t,hv)-\frac h2|v|^2
 =-\frac h2G_\lambda(v),\qquad
 G_\lambda(v)=\sup_{n,\pm}\{|v-p_{n,\pm}|^2-2\lambda b_n\}.
\]
This function is coercive and strongly convex, so its minimizer is unique and continuous in $\lambda>0$. If $v=(r,w)$, the maximum over the two signs is
\[
 r^2+2c_n|r|+c_n^2+(w-a_n)^2-2\lambda b_n.
\]
It is strictly larger at every $r\ne0$ than at $r=0$. Hence the minimizer is $(0,w(\lambda))$, where
\[
 w(\lambda)=\argmin_w\sup_n R_n(w,\lambda),\qquad
 R_n(w,\lambda)=c_n^2+(w-a_n)^2-2\lambda b_n.
\]

The intersection of consecutive parabolas occurs at
\begin{equation}\label{eq:counter-intersection}
 z_n(\lambda)=\frac{a_n+a_{n+1}}2+\frac{\lambda-\lambda_n}{e_n},\qquad
 R_n-R_{n+1}=2d_n(w-z_n(\lambda)).
\end{equation}
Set $\rho_n=e_nd_n/2=(A/2)n(n+1)4^{-n}$ and $J_n=[\lambda_n-\rho_n,\lambda_n+\rho_n]$. On $J_n$, $z_n$ runs from $a_n$ to $a_{n+1}$. For $n\ge5$, we have $e_k\ge e_n$ for $k<n$ and $e_k\le e_{n+1}<e_n-\rho_n$ for $k>n$. Thus for $\lambda\in J_n$,
\[
 z_k(\lambda)\ge\frac{a_k+a_{k+1}}2+1-\frac{d_n}2>1\quad(k<n),\qquad
 z_k(\lambda)\le\frac{a_k+a_{k+1}}2-\frac{e_n-\rho_n}{e_{n+1}}<0\quad(k>n).
\]
Equation~\eqref{eq:counter-intersection} therefore shows, at $w=z_n(\lambda)$,
\[
 R_1<\cdots<R_n=R_{n+1}>R_{n+2}>\cdots,
\]
including strict domination of the limiting parabola. Since $z_n\in[a_n,a_{n+1}]$, the derivatives of the two active parabolas are respectively nonnegative and nonpositive. Their convex hull contains zero, so strong convexity implies
\[
 w(\lambda)=z_n(\lambda)\quad(\lambda\in J_n).
\]
For $\lambda_n^\pm=\lambda_n\pm\rho_n\to\lambda_*$ and every $1/2<\alpha\le1$, we thus have
\begin{equation}\label{eq:counter-blowup}
 \frac{|w(\lambda_n^+)-w(\lambda_n^-)|}{|\lambda_n^+-\lambda_n^-|^\alpha}
 =\frac{2^{-n}}{(An(n+1)4^{-n})^\alpha}
 =A^{-\alpha}[n(n+1)]^{-\alpha}2^{(2\alpha-1)n}\longrightarrow\infty.
\end{equation}
The curve in \eqref{eq:counter-curve} is
$y(t)=(0,w(\lambda)/\lambda)$, $\lambda=(t-s)^{-1}$.
Near $\lambda_*>0$, the change of parameter and multiplication by $1/\lambda$ preserve local $\alpha$-H\"older regularity in both directions. Thus \eqref{eq:counter-blowup} disproves local $\alpha$-H\"older continuity at $t_*=s+1/(9A)$ for every $\alpha>1/2$.

To obtain a torus example, choose a small Euclidean ball in a flat coordinate chart and extend $\phi$ from a smaller ball to a Lipschitz function on the torus with a uniform Lipschitz constant $\ell$. Such an extension is given by the usual infimum extension formula. The slopes above are bounded independently of $A,s$. For the quadratic Hopf--Lax formula any minimizing initial point is within $2\ell t$ of its endpoint, and any maximizer in \eqref{eq:counter-curve} is within $2\ell(t-s)$ of the origin. First choose $s$ small, then $A$ large with $A^{-1}\ll s$. All points and optimizing segments for times near $t_*$ then remain in the chart, where \eqref{eq:counter-solution} is unchanged. Uniform Lipschitz and positive-time semiconcavity bounds give a short-time regularization threshold independent of $A$, so taking $A$ large also puts $t_*-s$ in this interval. This proves that local Lipschitz continuity can fail and, more strongly, that the local $1/2$-H\"older exponent in Theorem~\ref{thm:ISC} is optimal, even on a closed flat manifold.
\end{Rem}

\begin{Rem}\label{rem:frozen}
For an autonomous Hamiltonian and a fixed semiconcave profile $\phi$, the corresponding short-time curve
\[
 y^x(t)=\argmax_y\{\phi(y)-\breve h(t,x,y,\phi(y))\},\qquad y^x(0)=x,
\]
is Lipschitz. Here $\breve h(t,x,y,u)=\breve h(0,t,x,y,u)$. To see this, take smooth approximations $\phi_m\to\phi$ with uniform Lipschitz and semiconcavity bounds, and put
$\Psi_m(t,y)=\phi_m(y)-\breve h(t,x,y,\phi_m(y))$.
On the localization cone,
\[
 D_{yy}^2\Psi_m\le-\frac ctI,\qquad
 |\partial_tD_y\Psi_m|\le\frac Ct.
\]
The first estimate is Lemma~\ref{lem:strictconcavity}. The second follows from the endpoint formulas and their differentiated short-time estimates. It uses only $D\phi_m$, not a bound on its negative Hessian. The implicit function theorem at the unique maximizing point gives
\[
 \dot y_m^x(t)=-[D^2_{yy}\Psi_m(t,y_m^x(t))]^{-1}
                         \partial_tD_y\Psi_m(t,y_m^x(t)),
 \qquad |\dot y_m^x(t)|\le C.
\]
Localization extends this bound to $t=0$. Uniform convergence and uniqueness identify the limit with $y^x$, which is therefore Lipschitz. For the moving profile $u(t,\cdot)$ an additional mixed derivative of the profile occurs, and this proof does not apply. Remark~\ref{rem:nonLip} shows that the moving-profile curve need not be Lipschitz.
\end{Rem}

\begin{The}\label{thm:propagation}
If $(s,x)\in\Cut(u)$, then the intrinsic curve \eqref{eq:ISC} satisfies
\[
 (\tau,y^x(\tau))\in\Sing(u)\qquad\text{for every }s<\tau\le t.
\]
\end{The}
\begin{proof}
Suppose $u$ were differentiable at $(\tau,y^x(\tau))$. The positive optimizing lift in Theorem~\ref{thm:ISC} would have terminal momentum $D_xu(\tau,y^x(\tau))$, by the endpoint relation in the proof of Theorem~\ref{thm:Arnaud}. An incoming calibrated lift at the same endpoint has that same momentum and scalar value. Uniqueness of the contact flow identifies the two lifts on $[s,\tau]$. This provides a forward calibrated extension from $(s,x)$, contrary to the definition of the cut locus.
\end{proof}

\begin{Cor}\label{cor:Cutclosure}
For $u=T_0^t\phi$,
$\Cut(u)\subset\overline{\Sing(u)}$, where the closure is taken in $(0,\infty)\times M$.
\end{Cor}
\begin{proof}
The singular points $(\tau,y^x(\tau))$ converge to $(s,x)$ as $\tau\downarrow s$.
\end{proof}

\subsection{Lax--Oleinik commutators}\label{subsec:commutators}
For this subsection only, it is convenient to retain the auxiliary notation
\begin{equation}\label{eq:Aaux}
 A(a,b,x,y,u):=u+h(a,b,x,y,u),\qquad
 \breve A(a,b,x,y,u):=u-\breve h(a,b,x,y,u).
\end{equation}
These are scalar endpoint values. The fundamental solutions remain $h$ and $\breve h$.

\begin{Lem}\label{lem:inverse}
For fixed $a<b,x,y$, the scalar maps $A$ and $\breve A$ are mutually inverse increasing homeomorphisms. In particular, the endpoint inequalities
\[
 A(a,b,x,y,\breve A(a,b,x,y,u))\ge u,\qquad
 \breve A(a,b,x,y,A(a,b,x,y,u))\le u
\]
hold, and in fact both are equalities.
\end{Lem}
\begin{proof}
For an admissible curve $\xi$, write $E_\xi(r)$ for the terminal value of its scalar equation with initial value $r$. Scalar comparison gives
\begin{equation}\label{eq:Eorder}
 e^{-K(b-a)}(r_2-r_1)\le E_\xi(r_2)-E_\xi(r_1)
 \le e^{K(b-a)}(r_2-r_1),\quad r_1\le r_2.
\end{equation}
Thus $E_\xi$ is an increasing homeomorphism, and
$A(r)=\inf_\xi E_\xi(r)$, $\breve A(r)=\sup_\xi E_\xi^{-1}(r)$.
The two-sided estimate passes to $A$, which is consequently onto. Put $r_0=A^{-1}(r)$. For every $\xi$ we have $E_\xi(r_0)\ge r$, hence $E_\xi^{-1}(r)\le r_0$. Equality holds for a minimizer of $A(r_0)$. Therefore $\breve A(r)=r_0$. This also shows that the two endpoint problems have the same optimizing curves at inverse scalar values.
\end{proof}

\begin{Pro}\label{prop:commutators}
For $\phi\in C(M)$ and $a<b$,
\begin{enumerate}
\item $T_a^b\breve T_a^b\phi\ge\phi$ and $\breve T_a^bT_a^b\phi\le\phi$.
\item $T_a^b\breve T_a^bT_a^b\phi=T_a^b\phi$ and
$\breve T_a^bT_a^b\breve T_a^b\phi=\breve T_a^b\phi$.
\end{enumerate}
\end{Pro}
\begin{proof}
Endpoint inversion gives, for continuous $f,g$,
\[
 \breve A(a,b,x,y,f(y))\le g(x)
 \quad\Longleftrightarrow\quad
 f(y)\le A(a,b,x,y,g(x)).
\]
Quantifying over $x,y$ yields
$\breve T_a^bf\le g$ if and only if $f\le T_a^bg$.
Taking $g=\breve T_a^bf$, or $f=T_a^bg$, proves the two inequalities. Applying order preservation to these inequalities gives both directions of each triple-composition identity.
\end{proof}

\begin{Pro}\label{prop:reachability}
For $\phi\in C(M)$, the following conditions are equivalent.
\begin{enumerate}
\item $T_a^b\breve T_a^b\phi=\phi$.
\item there exists $\psi\in C(M)$ such that $T_a^b\psi=\phi$.
\item for every $x\in M$ there exists $y\in M$ such that
\[
 \breve T_a^b\phi(y)=\phi(x)-\breve h(a,b,y,x,\phi(x)).
\]
\end{enumerate}
\end{Pro}
\begin{proof}
The implication (i)$\Rightarrow$(ii) uses $\psi=\breve T_a^b\phi$, and (ii)$\Rightarrow$(i) is the first triple identity. For (i)$\Leftrightarrow$(iii), take an attained minimizing initial point in $T_a^b\breve T_a^b\phi(x)$ and apply Lemma~\ref{lem:inverse}. Conversely, the equality in (iii), followed by endpoint inversion, gives a competitor attaining $\phi(x)$. Proposition~\ref{prop:commutators} gives the opposite inequality.
\end{proof}

\begin{Rem}\label{rem:dual-equality}
The dual pointwise equality has an analogous form.
\[
 \breve T_a^bT_a^b\phi(x)=\phi(x)
 \quad\Longleftrightarrow\quad
 \exists y\in M:\ T_a^b\phi(y)=\phi(x)+h(a,b,x,y,\phi(x)).
\]
Indeed, use an attained positive maximizing endpoint and the same inverse relation. This is also the operator form of the calibration criterion in Lemma~\ref{lem:cal-existence}.
\end{Rem}

\begin{Pro}\label{prop:continuation}
Let $\phi$ be semiconcave and suppose $\breve T_a^b\phi\in C^{1,1}(M)$. Then the following statements hold.
\begin{enumerate}
\item $\breve T_t^b\phi\in C^{1,1}(M)$ for every $a\le t<b$.
\item $\Phi_H^{t,b}(\graph(D\breve T_t^b\phi))=\graph(D^+\phi)$ for every $a\le t<b$.
\item $\breve T_t^b\phi=T_a^t\breve T_a^b\phi$ and
$T_t^b\breve T_t^b\phi=\phi$ for every $a\le t\le b$.
\end{enumerate}
\end{Pro}
\begin{proof}
Set $g=\breve T_a^b\phi$. Fix $a<t\le b$ and a negative optimizer for $T_a^tg(z)$, with initial point $x$. First variation at its free initial point gives initial momentum $Dg(x)$. Choose a positive optimizer for $\breve T_a^b\phi(x)$. Its initial momentum is also $Dg(x)$. The latter assertion does not require smoothness of the long-time fundamental solution, because varying the initial endpoint of the optimizing curve gives a smooth lower support for $g$, whose derivative is the initial momentum. Differentiability of $g$ identifies the two.

Both lifts thus start at $(x,Dg(x),g(x))$. Uniqueness of \eqref{eq:contact} identifies the negative optimizer with the first part of the positive one. Dynamic programming along that curve gives
\[
 T_a^tg(z)=\breve T_t^b\phi(z).
\]
At $t=b$ this proves $T_a^bg=\phi$. For $a<t<b$ the common function is semiconcave as a negative evolution and semiconvex as a positive evolution, hence $C^{1,1}$. The remaining operator identity follows from the Markov property.

The same endpoint argument, now between $a\le t_1<t_2<b$, gives
\[
 \Phi_H^{t_1,t_2}(\graph(D\breve T_{t_1}^b\phi))
 =\graph(D\breve T_{t_2}^b\phi).
\]
One inclusion follows along each positive optimizer. The reverse follows by taking a negative optimizer ending at any point at $t_2$. Choose $t_2$ sufficiently close to $b$ for Theorem~\ref{thm:Arnaud}, applied to the fixed semiconcave profile $\phi$. Composing the flow identities proves (ii) for every $t_1$, including $a$.
\end{proof}

\section{Cut locus of contact systems}\label{sec:cut}
This section relates the geometry of the cut-time function to reachability
under the backward contact Hamiltonian flow, and then uses this viewpoint
to study the controllability of prescribed terminal data.

\subsection{Topological structure of the cut locus}
Continue to write $u(t,x)=T_0^t\phi(x)$. Define the cut time by
\begin{equation}\label{eq:cuttime}
 \tau_u(t,x):=\sup\{\tau\ge0:\text{a calibrated curve on }[t,t+\tau]
                                          \text{ starts at }x\}.
\end{equation}
Thus $\Cut(u)=\{(t,x):t>0,\ \tau_u(t,x)=0\}$. The value $\tau_u(t,x)$
is a remaining calibration time and may be infinite. If
$0<\tau_u(t,x)<+\infty$, then the supremum in \eqref{eq:cuttime} is attained.
Indeed, taking calibrated curves whose lengths increase to $\tau_u(t,x)$,
the uniform a priori estimates for minimizing curves and the Ascoli--Arzel\`a theorem yield, up to a subsequence, a limiting calibrated curve on
$[t,t+\tau_u(t,x)]$. When $t=0$, the same argument applies since the initial
point $x$ and the initial scalar value $u(0,x)=\phi(x)$ are fixed. Hence no
continuity of $\phi$ is required.

For a fixed $b>0$, set
\begin{equation}\label{eq:vb}
 v_b(t,x):=\breve T_t^b[u(b,\cdot)](x),\qquad 0\le t\le b.
\end{equation}
Domination gives $v_b\le u$. Recall $\mathcal R_u=\Pi_xD^*u$ from \eqref{eq:Ru}, and define
\begin{align*}
 G_b^*(u)&:=\{(y,p,u(b,y)):y\in M,\ p\in\mathcal R_u(b,y)\},\\
 G_b^\#(u)&:=\{(y,p,u(b,y)):y\in M,\ p\in D_y^+u(b,y)\setminus\mathcal R_u(b,y)\},\\
 \mathcal F_b^*(u)&:=\{\gamma:[0,b]\to M:\gamma(s)=\pi_x\Phi_H^{b,s}(y,p,u(b,y)),\
                                      (y,p,u(b,y))\in G_b^*(u)\},\\
 \mathcal A_b^*(u)&:=\{(t,x):0\le t\le b,\ \gamma(t)=x
                                      \text{ for some }\gamma\in\mathcal F_b^*(u)\},\\
 \mathcal A_b^\#(u)&:=([0,b]\times M)\setminus\mathcal A_b^*(u).
\end{align*}
Theorem~\ref{thm:cal-lifts} identifies $\mathcal F_b^*(u)$ with all calibrated curves on $[0,b]$. These definitions use $\mathcal R_u$, not $D_x^*u$. For example, for $H(p)=p^2/2$ the local solution
\[
 u(t,x)=\min\{0,-x-(t-b)/2,x-(t-b)/2\}
\]
has $\mathcal R_u(b,0)=\{-1,0,1\}$, whereas $D_x^*u(b,0)=\{-1,1\}$. The zero momentum is the calibrated branch arriving from the flat region at times less than $b$. Thus spatial reachable gradients alone need not parametrize all incoming calibrated curves.

\begin{Pro}\label{prop:contact-curves}
For the preceding sets and $v_b$, the following hold.
\begin{enumerate}
\item If $(t,x)\in\mathcal A_b^*(u)$ and $0<t<b$, there is a unique $\gamma\in\mathcal F_b^*(u)$ with $\gamma(t)=x$. Along its lift,
\begin{align*}
 &v_b(t,x)=u(t,x)=u(b,\gamma(b))-\int_t^bL(s,\gamma(s),\dot\gamma(s),u(s,\gamma(s)))\,ds,\\
 &D_xv_b(t,x)=D_xu(t,x)=p(t),\\
 &D_tv_b(t,x)+H(t,x,D_xv_b(t,x),v_b(t,x))=0.
\end{align*}
\item If $(t,x)\in\mathcal A_b^\#(u)$, every positive optimizer $\eta:[t,b]\to M$ for $v_b(t,x)$, with dual arc $p$, satisfies
\[
 (\eta(b),p(b),u(b,\eta(b)))\in G_b^\#(u),\qquad
 \eta(b)\in\Sing(u(b,\cdot)),\qquad v_b(t,x)<u(t,x).
\]
\end{enumerate}
\end{Pro}
\begin{proof}
Along a curve in $\mathcal F_b^*(u)$, calibration and domination give $v_b=u$ and the displayed integral formula. At an interior point, $u$ is semiconcave and $v_b$ semiconvex. Proposition~\ref{prop:touch} gives equality of their full gradients. Theorem~\ref{thm:cal-lifts} identifies the spatial gradient with $p(t)$, and the equation determines the time derivative. The fixed state $(x,D_xu(t,x),u(t,x))$ determines a unique contact lift, proving uniqueness of the curve.

We justify the endpoint relation for a positive optimizer also when $b-t$ is not small. Put $y=\eta(b)$ and choose $c<b$ sufficiently close to $b$. Keeping $\eta$ fixed on $[t,c]$, vary only its terminal point on the short interval $[c,b]$. By the short-time regularity of the fundamental solution, the resulting terminal scalar value is a smooth upper support for $u(b,\cdot)$ at $y$. Otherwise its inverse scalar evolution would give a larger value in the supremum defining $v_b(t,x)$. The first variation formula then gives $p(b)\in D_y^+u(b,y)$.

If $p(b)\in\mathcal R_u(b,y)$, Theorem~\ref{thm:cal-lifts} and flow uniqueness would make $\eta$ part of an incoming calibrated curve, contradicting $(t,x)\in\mathcal A_b^\#(u)$. Hence its terminal state lies in $G_b^\#(u)$. If $u(b,\cdot)$ were differentiable at $y$, its superdifferential would be a singleton containing the nonempty set $\mathcal R_u(b,y)$, again a contradiction. Finally, if $v_b(t,x)=u(t,x)$, then $\eta$ is calibrated by Lemma~\ref{lem:cal-existence}. For $t>0$, concatenating it with an incoming calibrated curve gives a calibrated curve on $[0,b]$. For $t=0$, $\eta$ itself already has this property. Hence $(t,x)\in\mathcal A_b^*(u)$, a contradiction. Thus $v_b(t,x)<u(t,x)$.
\end{proof}

\begin{Pro}\label{prop:contact-cuttime}
For every $0\le t\le b$ and $x\in M$,
\begin{align}
 (t,x)\in\mathcal A_b^*(u)
 &\Longleftrightarrow v_b(t,x)=u(t,x)
 \Longleftrightarrow\tau_u(t,x)\ge b-t,\label{eq:Astar}\\
 (t,x)\in\mathcal A_b^\#(u)
 &\Longleftrightarrow v_b(t,x)<u(t,x)
 \Longleftrightarrow\tau_u(t,x)<b-t.\label{eq:Ahash}
\end{align}
Moreover, $Dv_b=Du$ is locally Lipschitz on
$\mathcal A_b^*(u)\cap((0,b)\times M)$.
\end{Pro}
\begin{proof}
A calibrated segment from $(t,x)$ to time $b$ gives equality in \eqref{eq:vb}. Conversely, an attained positive optimizer at an equality point is calibrated. For $t>0$, concatenate it with an incoming calibrated curve. For $t=0$, it already has the required full interval. This proves the first equivalence. Attainment of finite positive cut times proves the second, and $t=b$ is immediate. The complementary equivalences use $v_b\le u$. The gradient assertion follows from Proposition~\ref{prop:touch} on compact interior space-time charts.
\end{proof}

Define the regularization length restricted to nonnegative initial times by
\[
 \bar\rho_1(u(b,\cdot)):=\sup\{b-a:0\le a<b,\ \breve T_a^bu(b,\cdot)\in C^{1,1}(M)\}.
\]
It is positive by Theorem~\ref{thm:Arnaud}. Proposition~\ref{prop:continuation} ensures regularity at every shorter positive length.

\begin{The}\label{thm:topology}
For every $0<t<b$, the sets
$\{x\in M:\tau_u(t,x)\ge b-t\}$ and $G_b^*(u)$ are bi-Lipschitz homeomorphic. The maps are
\begin{equation}\label{eq:topology-maps}
 (y,p,u(b,y))\longmapsto\pi_x\Phi_H^{b,t}(y,p,u(b,y)),\qquad
 x\longmapsto\Phi_H^{t,b}(x,D_xu(t,x),u(t,x)).
\end{equation}
If $b-\bar\rho_1(u(b,\cdot))<t<b$, the same backward projection is a bi-Lipschitz homeomorphism from $G_b^\#(u)$ onto
$\{x\in M:\tau_u(t,x)<b-t\}$, with inverse
\[
 x\longmapsto\Phi_H^{t,b}(x,D_xv_b(t,x),v_b(t,x)).
\]
\end{The}
\begin{proof}
Propositions~\ref{prop:contact-curves} and \ref{prop:contact-cuttime} give bijectivity and the inverse formula in \eqref{eq:topology-maps}. The set $G_b^*(u)$ is compact and its trajectories are calibrated on $[0,b]$. The contact flow is $C^1$ on a neighborhood of this compact family, so the backward projection is Lipschitz. On the compact contact slice $\{u(t,\cdot)=v_b(t,\cdot)\}$, Proposition~\ref{prop:touch} gives a Lipschitz common gradient. The second map in \eqref{eq:topology-maps} is therefore Lipschitz as well.

Under the additional time restriction, Proposition~\ref{prop:continuation} gives
\[
 \Phi_H^{b,t}(\graph(D_y^+u(b,\cdot)))=\graph(D_xv_b(t,\cdot)).
\]
The first part identifies the image of $G_b^*(u)$ with the portion over the contact slice. Injectivity of the flow identifies the complementary portions. Both directions are Lipschitz because $v_b(t,\cdot)\in C^{1,1}(M)$ and the flow has uniform derivative bounds on the relevant compact set.
\end{proof}

\begin{Rem}\label{rem:dimension}
The maps in Theorem~\ref{thm:topology} preserve Hausdorff dimension and countable $k$-rectifiability for every $k$. In particular, all the interior slices $\{\tau_u(t,\cdot)\ge b-t\}$ are mutually bi-Lipschitz. The same is true of the complementary slices when the times lie in the regularization interval. For comparison, the singular strata of a semiconcave time slice,
\[
 \Sigma_k(u(t,\cdot)):=\{x\in M:\dim D_x^+u(t,x)=k\},
 \qquad k=1,\ldots,m,
\]
are countably $(m-k)$-rectifiable (see \cite{Cannarsa_Sinestrari_book}). Together with
$\Sing(u)\subset\Cut(u)\subset\overline{\Sing(u)}$, this relates the cut-time problem to the geometry of terminal momentum sets. The bi-Lipschitz correspondences do not alone imply rectifiability of the entire cut locus. That would require further information relating $G_b^*(u)$ and $G_b^\#(u)$ to these singular strata.
\end{Rem}

\subsection{Controllability of contact systems}\label{subsec:control}
Let $\phi$ now be a prescribed semiconcave terminal profile at time $b$. Define
\begin{align}
 \rho_1^b(\phi)&:=\sup\{r>0:\breve T_{b-r}^b\phi\in C^{1,1}(M)\},\label{eq:rho1}\\
 \rho_2^b(\phi)&:=\sup\{r\ge0:T_{b-r}^b\breve T_{b-r}^b\phi=\phi\}.\label{eq:rho2}
\end{align}
The Hamiltonian is defined for all real times, so $b-r$ is allowed to be negative here. One has $0<\rho_1^b(\phi)\le\rho_2^b(\phi)$. Regularization holds for every $0<r<\rho_1^b(\phi)$. No endpoint regularity is required. Reachability holds exactly for $r\in[0,\rho_2^b(\phi)]$ when $\rho_2^b(\phi)$ is finite, and for all $r\ge0$ otherwise. Indeed, Proposition~\ref{prop:continuation} proves persistence of regularity under shortening and shows that every regularization length is also a reachable length. The Markov property preserves reachability under shortening, while continuity of the operators in their time endpoints makes the set of finite reachable lengths closed.

Fix $a<b$ with $b-a\le\rho_2^b(\phi)$, so $T_a^b\breve T_a^b\phi=\phi$. Set
\begin{equation}\label{eq:ubar-breve}
 \breve u(t,x):=\breve T_t^b\phi(x),\qquad
 \bar u(t,x):=T_a^t\breve T_a^b\phi(x),\qquad a\le t\le b.
\end{equation}
Define $\mathcal F(\phi;a,b)$ to be the family of $C^1$ curves $\gamma:[a,b]\to M$ satisfying
\begin{equation}\label{eq:Fcontrol}
 \phi(\gamma(b))=\breve T_a^b\phi(\gamma(a))
             +\int_a^bL(s,\gamma(s),\dot\gamma(s),u_\gamma(s))\,ds,
 \qquad u_\gamma(a)=\breve T_a^b\phi(\gamma(a)).
\end{equation}
Let $\mathcal A(\phi;a,b)$ be their union in space-time. Also put
\begin{align*}
 \mathcal F^*(\phi;a,b)&:=\{\gamma(s)=\pi_x\Phi_H^{b,s}(x,p,\phi(x)):
                            x\in M,\ p\in D^*\phi(x),\ a\le s\le b\},\\
 \mathcal A^*(\phi;a,b)&:=\{(t,x):\gamma(t)=x
                              \text{ for some }\gamma\in\mathcal F^*(\phi;a,b)\}.
\end{align*}
Here $D^*\phi$ is the spatial reachable-gradient set of the prescribed profile. Its role differs from $\mathcal R_u$ in the preceding subsection.

\begin{Lem}\label{lem:universal}
If $\psi\in C(M)$ satisfies $T_a^b\psi=\phi$, every curve in $\mathcal F^*(\phi;a,b)$ is calibrated for $u_\psi(t,x):=T_a^t\psi(x)$.
\end{Lem}
\begin{proof}
For $p\in D^*\phi(x)$ choose spatial differentiability points $x_j\to x$ with $D\phi(x_j)\to p$. Each incoming $u_\psi$-calibrated lift at $(b,x_j)$ has terminal momentum $D\phi(x_j)$, by the terminal supergradient relation. The minimizing-curve bounds and continuous dependence give the lift issued backward from $(x,p,\phi(x))$. Continuity of $u_\psi$ on $[a,b]\times M$ passes calibration to the limit. This proves at the same time that each such backward trajectory exists on the whole interval.
\end{proof}

\begin{Pro}\label{prop:min-control}
The following hold for the minimal initial datum $\breve T_a^b\phi$.
\begin{enumerate}
\item Every $(t,x)\in\mathcal A(\phi;a,b)$ with $a<t<b$ lies on a unique curve $\gamma\in\mathcal F(\phi;a,b)$. Along its lift,
\begin{align*}
 \bar u(t,x)=\breve u(t,x)
 &=\phi(\gamma(b))-\int_t^bL(s,\gamma,\dot\gamma,u_\gamma)\,ds,\\
 D_x\bar u(t,x)&=D_x\breve u(t,x)=p(t).
\end{align*}
\item $\mathcal A(\phi;a,b)=\{(t,x):\bar u(t,x)=\breve u(t,x)\}$. The common full gradient is locally Lipschitz on its interior part. If $b-a<\rho_1^b(\phi)$, then $\mathcal A(\phi;a,b)=[a,b]\times M$.
\item $\mathcal F^*(\phi;a,b)\subset\mathcal F(\phi;a,b)$ and
$\mathcal A^*(\phi;a,b)\subset\mathcal A(\phi;a,b)$.
\end{enumerate}
\end{Pro}
\begin{proof}
Domination gives $\bar u\ge\breve u$. Every curve satisfying \eqref{eq:Fcontrol} is $\bar u$-calibrated, and comparison from the two endpoints gives equality $\bar u=\breve u$ along it and the integral identity in (i).

Conversely, suppose $\bar u(t,x)=\breve u(t,x)$ with $a<t<b$. Concatenate a negative optimizer for $\bar u(t,x)$ and a positive optimizer for $\breve u(t,x)$. Their scalar values agree at the junction. The concatenated curve starts from the minimal datum and ends at the value $\phi$, hence is minimizing for the full negative evolution. Herglotz regularity makes it a single $C^2$ calibrated curve. At $t=a$ one uses a positive optimizer, and at $t=b$ a negative optimizer. This proves the equality-set assertion.

On the interior equality set, semiconcavity of $\bar u$ and semiconvexity of $\breve u$ imply equality and local Lipschitz continuity of their full gradients by Proposition~\ref{prop:touch}. Calibration identifies the spatial gradient with $p(t)$. Uniqueness of the contact flow then fixes the whole curve through an interior point. Proposition~\ref{prop:continuation} gives $\bar u=\breve u$ everywhere when $b-a<\rho_1^b(\phi)$. Finally, (iii) follows from Lemma~\ref{lem:universal}, applied to $\psi=\breve T_a^b\phi$.
\end{proof}

For any other continuous $\psi$ reaching $\phi$, let $\mathcal F(\phi,\psi;a,b)$ be defined by \eqref{eq:Fcontrol} with initial value $\psi(\gamma(a))$ in place of $\breve T_a^b\phi(\gamma(a))$, and let $\mathcal A(\phi,\psi;a,b)$ be its space-time union.

\begin{Pro}\label{prop:all-control-sets}
If $T_a^b\psi=\phi$, then
\begin{align*}
 \mathcal F^*(\phi;a,b)&\subset\mathcal F(\phi,\psi;a,b)\subset\mathcal F(\phi;a,b),\\
 \mathcal A^*(\phi;a,b)&\subset\mathcal A(\phi,\psi;a,b)\subset\mathcal A(\phi;a,b).
\end{align*}
Moreover $u_\psi\ge\bar u\ge\breve u$, and
\[
 \mathcal A(\phi,\psi;a,b)
 =\{(t,x):u_\psi(t,x)=\bar u(t,x)=\breve u(t,x)\}.
\]
The common full gradient is locally Lipschitz on $\mathcal A(\phi,\psi;a,b)\cap((a,b)\times M)$. More precisely,
\begin{equation}\label{eq:familycontrol}
 \mathcal F(\phi,\psi;a,b)
 =\{\gamma\in\mathcal F(\phi;a,b):
                   \psi(\gamma(a))=\breve T_a^b\phi(\gamma(a))\}.
\end{equation}
\end{Pro}
\begin{proof}
Proposition~\ref{prop:commutators} gives
$\psi\ge\breve T_a^bT_a^b\psi=\breve T_a^b\phi$. Order preservation implies $u_\psi\ge\bar u\ge\breve u$. Along $\gamma\in\mathcal F(\phi,\psi;a,b)$, invert its scalar equation to obtain
\[
 \psi(\gamma(a))=E_\gamma^{-1}(\phi(\gamma(b)))
 \le\breve T_a^b\phi(\gamma(a))\le\psi(\gamma(a)).
\]
Thus the initial values agree and $\gamma\in\mathcal F(\phi;a,b)$. The converse in \eqref{eq:familycontrol} is immediate. Lemma~\ref{lem:universal} gives the first family inclusion. The set inclusions follow.

On each such curve, comparison from the two endpoints gives equality of all three functions. At any interior point of triple equality, concatenate a negative optimizer for $u_\psi$ and a positive optimizer for $\breve u$. As in Proposition~\ref{prop:min-control}, the concatenation is minimizing for the full evolution and belongs to $\mathcal F(\phi,\psi;a,b)$. The endpoint cases follow by taking just the positive or negative optimizer. Finally apply Proposition~\ref{prop:touch} to $\breve u\le u_\psi$ and to $\breve u\le\bar u$ to obtain the common full gradient and its local Lipschitz bound on the interior-time equality set.
\end{proof}

\begin{Cor}\label{cor:controls}
For $\psi\in C(M)$, the following conditions are equivalent.
\begin{enumerate}
\item $T_a^b\psi=\phi$.
\item the forward viscosity solution with initial datum $\psi$ solves
\[
 \partial_tu_\psi+H(t,x,D_xu_\psi,u_\psi)=0,
 \qquad u_\psi(a,\cdot)=\psi,\quad u_\psi(b,\cdot)=\phi.
\]
\item $\psi\ge\breve T_a^b\phi$ on $M$, with equality on
\[
 \pi_x\bigl(\mathcal A^*(\phi;a,b)\cap(\{a\}\times M)\bigr).
\]
\end{enumerate}
\end{Cor}
\begin{proof}
The equivalence of (i) and (ii) follows from the variational representation and uniqueness of the forward Cauchy problem. If (i) holds, the commutator inequality gives $\psi\ge\breve T_a^b\phi$. Lemma~\ref{lem:universal} and \eqref{eq:familycontrol} force equality at the initial point of every curve in $\mathcal F^*(\phi;a,b)$, proving (iii).

Conversely, (iii) gives $T_a^b\psi\ge T_a^b\breve T_a^b\phi=\phi$. At an arbitrary terminal point $x$, choose $p\in D^*\phi(x)$ and the backward lift $\gamma$ issued from $(x,p,\phi(x))$. By Proposition~\ref{prop:min-control}(3), $\gamma\in\mathcal F(\phi;a,b)$. Its initial point lies in the set in (iii), so $\psi(\gamma(a))=\breve T_a^b\phi(\gamma(a))$. Hence the same lift is a competitor for $T_a^b\psi(x)$ with terminal value $\phi(x)$, giving $T_a^b\psi(x)\le\phi(x)$. Thus $T_a^b\psi=\phi$.
\end{proof}
\begin{Rem}
Equivalently, every continuous control is $\breve T_a^b\phi+\eta$, where $\eta\ge0$ is continuous and vanishes on the initial projection in (iii). That projection is compact. Hence the continuous control is unique exactly when this initial projection equals $M$. Spatial reachable gradients suffice here because they supply at least one universally calibrated lift at each terminal point. This does not identify them with all the incoming calibrated momenta of a given solution.
\end{Rem}

\appendix
\section{Properties of contact fundamental solutions}\label{app:A}
We collect the short-time facts used in Section~\ref{sec:arnaud}. All constants are uniform on a fixed compact time interval $I$. Calculations are made in coordinate balls. A finite atlas transfers them to $M$. We use the Herglotz a priori estimates from \cite[Theorem~2.4 and Proposition~A.3]{CCJWY2020}. For small $t-s$, with $d(x,y)\le\lambda(t-s)$ and bounded prescribed scalar endpoint, every minimizing lift stays in a fixed compact subset of $TM\times\R$, equivalently of $T^*M\times\R$. The estimates are also uniform for bounded endpoint values when the time length is bounded above and away from zero.

For the short-time uniformity, comparison with a constant-speed short geodesic bounds the action and the integral of a superlinear function of the speed by $C(t-s)$. The time-dependent Erdmann estimate then bounds the momentum. Along a regular extremal its differential form is
\[
 \frac d{d\tau}H(\tau,\xi,p,u_\xi)=H_t-H_uH=-L_t+L_uH.
\]
The growth condition (L4), scalar comparison, and Gronwall's inequality give the required uniform bound. The variational regularity argument in \cite{CCJWY2020} justifies this estimate before regularity of a minimizer is known. Positive endpoint problems are reduced to negative ones by the following reversal.

\begin{Lem}\label{lem:reversal}
For fixed $s<t$, define
$\breve L(\tau,x,v,u):=L(s+t-\tau,x,-v,-u)$.
Then
\begin{equation}\label{eq:reversal}
 \breve h_L(s,t,x,y,u)=h_{\breve L}(s,t,y,x,-u).
\end{equation}
\end{Lem}
\begin{proof}
If $\xi$ runs from $x$ to $y$ and $u_\xi(t)=u$, set
$\eta(\tau)=\xi(s+t-\tau)$ and $w(\tau)=-u_\xi(s+t-\tau)$.
Then $\eta$ runs from $y$ to $x$, $w(s)=-u$, and
$\dot w=\breve L(\tau,\eta,\dot\eta,w)$. Moreover,
$w(t)-w(s)=u_\xi(t)-u_\xi(s)$.
This bijection preserves the action, proving \eqref{eq:reversal}. The transformed Lagrangians obey the same bounds uniformly for $s,t\in I$.
\end{proof}

\begin{Pro}\label{prop:kernel}
Fix $\lambda,\delta>0$. There are $T_{\lambda,\delta}>0$ and $c,C>0$ such that for $0<t-s<T_{\lambda,\delta}$, $s,t\in I$, $|u|\le\delta$, and $d(x,y)\le\lambda(t-s)$, the following hold on slightly smaller cones.
\begin{enumerate}
\item Both $x\mapsto\breve h(s,t,x,y,u)$ and $y\mapsto\breve h(s,t,x,y,u)$ are locally semiconcave and uniformly convex, with Hessian bounds
\begin{equation}\label{eq:kernelHess}
 \frac c{t-s}I\le\breve h_{xx},\breve h_{yy}\le\frac C{t-s}I
\end{equation}
wherever these derivatives are evaluated. The minimizer $\xi$ is unique.
\item For fixed $s,x$, the function $(t,y,u)\mapsto\breve h(s,t,x,y,u)$ is locally $C^{1,1}$. Along its minimizing lift, with $p=L_v(\tau,\xi,\dot\xi,u_\xi)$ and
\[
 \mu=\exp\!\left(-\int_s^tL_u(\tau,\xi,\dot\xi,u_\xi)\,d\tau\right),
\]
one has
\begin{equation}\label{eq:terminalderivs}
 \breve h_t=-\mu H(t,y,p(t),u),\qquad
 \breve h_y=\mu p(t),\qquad \breve h_u=1-\mu.
\end{equation}
\item For fixed $t,y$, the function $(s,x,u)\mapsto\breve h(s,t,x,y,u)$ is locally $C^{1,1}$, and
\begin{equation}\label{eq:initialderivs}
 \breve h_s=H(s,x,p(s),u_\xi(s)),\qquad \breve h_x=-p(s).
\end{equation}
\end{enumerate}
\end{Pro}
\begin{proof}
The a priori bound places all terminal momenta of minimizers in a fixed ball. If two minimizers have terminal momenta $p_1,p_2$, their backward lifts have the same initial position $x$ and terminal data $(y,u)$. Lemma~\ref{lem:twist} forces $p_1=p_2$, proving uniqueness. Its inverse function argument gives a $C^1$ terminal momentum
$p=p(s,t,x,y,u)$ solving $X(s;t,y,p,u)=x$.

For a variation of a contact trajectory at a fixed intermediate time, direct differentiation of \eqref{eq:contact} gives
\begin{equation}\label{eq:contactvariation}
 \frac d{d\tau}(\delta U-P\cdot\delta X)
             =-H_u(\delta U-P\cdot\delta X).
\end{equation}
Put $B=u-\breve h$ temporarily. When the endpoint times vary, substitute
$\delta X(s)=dx-\dot X(s)\,ds$ and
$\delta U(s)=dB-\dot U(s)\,ds$, and the analogous expressions at $t$. Since $\dot U-P\cdot\dot X=-H$, integration gives
\begin{equation}\label{eq:fullfirstvariation}
 dB=P(s)\cdot dx-H(s)\,ds
                   +\mu\{du-P(t)\cdot dy+H(t)\,dt\}.
\end{equation}
This proves \eqref{eq:terminalderivs}--\eqref{eq:initialderivs}. In particular, the first derivatives are $C^1$, which gives the stated regularity.

Write $T=t-s$. The variational equation and the implicit equation for $p$ yield
\begin{align}
 p_y&=-X_p^{-1}X_y=T^{-1}H_{pp}(t,y,p,u)^{-1}+o(T^{-1}),\label{eq:py}\\
 p_x&=X_p^{-1}=-T^{-1}H_{pp}(t,y,p,u)^{-1}+o(T^{-1}),
 &p_u&=-X_p^{-1}X_u=O(1).\label{eq:pxpu}
\end{align}
All errors are uniform on the compact phase-space set. In particular, derivatives of the lift with respect to $x,y$ are $O(T^{-1})$ and with respect to $u$ are $O(1)$. If
$J=\int_s^t L_u\,d\tau=-\int_s^tH_u\,d\tau$, then
\begin{equation}\label{eq:Jbounds}
 J=O(T),\qquad J_x,J_y=O(1),\qquad J_u=O(T).
\end{equation}
Differentiating $\breve h_y=e^{-J}p(t)$ gives
$\breve h_{yy}=T^{-1}H_{pp}^{-1}+o(T^{-1})$.
Likewise $\breve h_x=-P(s)$, $P_p(s)=I+O(T)$, and \eqref{eq:pxpu} give the same leading term for $\breve h_{xx}$. Uniform positive definiteness proves \eqref{eq:kernelHess} and the local semiconcavity bounds. The same differentiation also gives $|\breve h_{xy}|\le C/T$.
\end{proof}

\begin{Rem}\label{rem:kernelC2}
The proof gives joint local $C^2$ regularity of $\breve h$ in $(s,t,x,y,u)$ on these short-time cones under the standing $C^2$ assumption on $H$. This follows because the terminal momentum is $C^1$ and the first-variation formulas have $C^1$ right sides. No third derivatives of $H$ are needed.

We record two additional consequences used above. First,
\begin{equation}\label{eq:mixedextra}
 |\breve h_{xu}|+|\breve h_{yu}|\le C,
 \qquad |\breve h_{uu}|\le CT,
 \qquad |\breve h_{ty}|+|\breve h_{tu}|\le C/T.
\end{equation}
Indeed, $p_s,p_t=O(T^{-1})$ by the implicit equation, while $J_s,J_t=O(1)$, including the bounded moving-endpoint terms. These and \eqref{eq:Jbounds} give \eqref{eq:mixedextra} directly. Second, the negative endpoint value $A=u+h$, now with $u=U(s)$, satisfies
\begin{equation}\label{eq:negativevariation}
 dA=P(t)\cdot dy-H(t)\,dt
       +\nu\{du-P(s)\cdot dx+H(s)\,ds\},\qquad
 \nu=\exp\!\left(\int_s^tL_u\,d\tau\right).
\end{equation}
It follows either from \eqref{eq:contactvariation} or from the scalar inverse relation. Reversal \eqref{eq:reversal} is useful for estimates, but the reversed Lagrangian depends on $s,t$. Ignoring that dependence when differentiating endpoint times would give an unjustified formula.
\end{Rem}

We next retain the smooth terminal construction. Let $H\in C^{r+1}$ and $\phi\in C^{r+1}(M)$, $r\ge1$, and write
\begin{equation}\label{eq:smooth-lift}
 (X(s,z),P(s,z),U(s,z))=\Phi_H^{t,s}(z,D\phi(z),\phi(z)).
\end{equation}
Equivalently, this is the terminal characteristic system
\begin{equation*}\label{eq:Liephi}\tag{Lie$_\phi$}
 \begin{aligned}
 \dot X&=H_p(s,X,P,U),&X(t)&=z,\\
 \dot P&=-H_x(s,X,P,U)-H_u(s,X,P,U)P,&P(t)&=D\phi(z),\\
 \dot U&=\langle P,H_p(s,X,P,U)\rangle-H(s,X,P,U),&U(t)&=\phi(z).
 \end{aligned}
\end{equation*}
The compact terminal graph has a backward flow on $[t-\sigma_1,t]$, uniformly when $\phi$ ranges over a bounded subset of $C^2(M)$ and $t$ ranges over $I$. Define $\Phi(s,z):=(s,X(s,z))$.

\begin{Lem}\label{lem:projection}
There is $\sigma_2=\sigma_2(\|\phi\|_{C^2})\in(0,\sigma_1]$, uniform for $t\in I$, such that the characteristic projection
\[
 \Phi:(s,z)\longmapsto(s,X(s,z))
\]
is a $C^r$ diffeomorphism of $[t-\sigma_2,t]\times M$ onto itself. Its inverse has the form $(s,x)\mapsto(s,Z(s,x))$.
\end{Lem}
\begin{proof}
The derivative $D_zX(s,z)$ is uniformly close to the identity when $t-s$ is small. Hence each spatial projection is a local diffeomorphism. Compactness makes it a covering map, and the homotopy to the identity implies that the induced map on the fundamental group is surjective. The covering therefore has one sheet, and the map is globally invertible. Injectivity also follows directly from the uniform local inverse estimates, which exclude nearby distinct preimages, and from the uniform closeness to the identity, which excludes distant ones. The inverse function theorem with parameter $s$ gives the asserted regularity and form of the inverse.
\end{proof}

\begin{The}\label{thm:classical}
For $\phi\in C^{r+1}(M)$ and $H\in C^{r+1}$, $r\ge1$, the terminal problem
\[
 w_s+H(s,x,D_xw,w)=0,\qquad w(t,x)=\phi(x),
\]
has on $[t-\sigma_2,t]\times M$ a unique classical solution $w\in C^{r+1}$, given by
\begin{equation}\label{eq:classicalrepresentation}
 w(s,x)=U(s,Z(s,x))=\breve T_s^t\phi(x),\qquad
 D_xw(s,X(s,z))=P(s,z).
\end{equation}
\end{The}
\begin{proof}
On the terminal graph, $d\phi-D\phi\cdot dz=0$. Thus \eqref{eq:contactvariation} implies $U_z=P X_z$ along \eqref{eq:smooth-lift}. Define $w(s,x)=U(s,Z(s,x))$. It follows that $D_xw=P(s,Z(s,x))$, and differentiating in time gives $w_s=-H(s,x,D_xw,w)$. These derivatives are $C^r$, so $w$ is $C^{r+1}$.

Along any admissible curve from time $s$ to $t$, Fenchel's inequality gives
$\frac d{d\tau}w(\tau,\xi(\tau))\le L(\tau,\xi,\dot\xi,w(\tau,\xi))$.
Scalar comparison from the terminal value yields $w(s,x)\ge\breve T_s^t\phi(x)$. The characteristic through $(s,x)$ attains equality, giving \eqref{eq:classicalrepresentation}. The same comparison also proves uniqueness among classical solutions.
\end{proof}

\begin{Pro}\label{prop:boundedC2}
For every bounded family $F\subset C^2(M)$ there is $t_F>0$, uniform in the time endpoints in $I$, such that $\breve T_s^t(F)$ is bounded in $C^2(M)$ whenever $0\le t-s\le t_F$. If $f\in F$ and
\[
 \eta(\tau)=\pi_x\Phi_H^{t,\tau}(y,Df(y),f(y)),
\]
then, along its scalar lift $w_\eta(t)=f(y)$,
\[
 \breve T_s^tf(\eta(s))=w_\eta(s)
 =f(y)-\int_s^tL(\tau,\eta(\tau),\dot\eta(\tau),w_\eta(\tau))\,d\tau.
\]
\end{Pro}
\begin{proof}
The compactness bounds for the terminal graphs, the variational equation, and the uniform inverse bound for $D_zX$ in Lemma~\ref{lem:projection} give a common $C^1$ bound on $P(s,Z(s,x))$. Theorem~\ref{thm:classical} therefore gives uniform $C^2$ bounds for $\breve T_s^tf$. Its representation along the characteristic and the scalar equation give the displayed identity.
\end{proof}

\begin{Lem}\label{lem:localization}
Let $\phi$ be Lipschitz with constant $\ell$, and let its norm be bounded. There is $\lambda_\ell>0$, uniform for these bounds and the fixed time interval $I$, such that every maximizing endpoint in $\breve T_s^t\phi(x)$ belongs to $B(x,\lambda_\ell(t-s))$ for $s<t$ in $I$.
\end{Lem}
\begin{proof}
Put $T=t-s\le |I|$. All constants in this proof may depend on $I$. A constant curve gives
$\breve T_s^t\phi(x)\ge\phi(x)-CT$.
For any admissible curve $\xi$ ending at $y$ with terminal scalar value $\phi(y)$, the scalar equation, superlinearity, and $|L_u|\le K$ give
\[
 u_\xi(s)\le\phi(y)+CT
       -e^{-KT}\int_s^t\theta_0(|\dot\xi|)\,d\tau.
\]
One can replace $\theta_0$ by a convex superlinear lower bound. At an optimizer, Jensen's inequality and $\phi(y)-\phi(x)\le\ell d(x,y)$ imply
\[
 e^{-KT}T\theta_0(d(x,y)/T)\le\ell d(x,y)+CT.
\]
Superlinearity bounds $d(x,y)/T$ uniformly, proving the claim. The constants also depend on the common scalar bound, which is fixed in every application.
\end{proof}

\begin{Lem}\label{lem:scalarmonotonicity}
Both $u\mapsto u+h(s,t,x,y,u)$ and
$u\mapsto u-\breve h(s,t,x,y,u)$ are increasing. More precisely, for either map $F$ and $u_1\le u_2$,
\[
 e^{-K(t-s)}(u_2-u_1)\le F(u_2)-F(u_1)
                         \le e^{K(t-s)}(u_2-u_1).
\]
\end{Lem}
\begin{proof}
For a fixed curve, the difference of two scalar solutions satisfies
$\dot w=a(\tau)w$ with $|a|\le K$. This proves the two-sided estimate both for its initial-to-terminal map and its inverse. Taking the infimum of the former or the supremum of the latter preserves the same two-sided estimate. In particular, the monotonicity used in all the operator arguments does not require a sign for $L_u$.
\end{proof}

\begin{Lem}\label{lem:mixed}
Fix a short-time cone and a bounded scalar interval. With
\[
 F_{s,t,x}(y,u):=u-\breve h(s,t,x,y,u),\qquad T=t-s,
\]
there are $C,T_0>0$, independent of $T\in(0,T_0]$, such that
\begin{align*}
 |D_yF_{s,t,x}(y,u_1)-D_yF_{s,t,x}(y,u_2)|&\le C|u_1-u_2|,\\
 |D_uF_{s,t,x}(y,u_1)-D_uF_{s,t,x}(y,u_2)|&\le CT|u_1-u_2|.
\end{align*}
\end{Lem}
\begin{proof}
The implicit terminal momentum equation gives
$p_u=-X_p^{-1}X_u=O(1)$, since $X_p^{-1}=O(T^{-1})$ and $X_u=O(T)$. The derivative of the entire lift with respect to the prescribed terminal scalar value is consequently bounded. Therefore, for $J=\int_s^t L_u\,d\tau$, one has $J_u=O(T)$.
The formulas
$F_u=e^{-J}$ and $F_y=-e^{-J}p(t)$ imply
$F_{uu}=O(T)$ and $F_{yu}=O(1)$.
Integrating these derivative bounds between $u_1$ and $u_2$ proves the assertion. This explicit factor $T$ in the second estimate is useful when composing with a merely semiconcave profile.
\end{proof}

\begin{Lem}\label{lem:strictconcavity}
Under the bounds of Lemma~\ref{lem:localization}, assume in addition that $\phi$ has semiconcavity constant $C_\phi$. For sufficiently small $T=t-s>0$, the function
\[
 y\longmapsto\phi(y)-\breve h(s,t,x,y,\phi(y))
\]
has semiconcavity constant at most $-c/T$ on a coordinate ball containing all its maximizing endpoints. In particular, its global maximizer is unique. The constants are uniform for uniformly bounded Lipschitz and semiconcave families.
\end{Lem}
\begin{proof}
Enlarge the cone in Lemma~\ref{lem:localization} slightly and write $F(y,u)=u-\breve h(s,t,x,y,u)$. Proposition~\ref{prop:kernel} and Lemma~\ref{lem:mixed} give
\[
 F_{yy}\le-\frac{c_0}T I,\quad |F_{yu}|\le C,\quad
 |F_{uu}|\le CT,\quad 0<F_u\le e^{KT}.
\]
Smooth $\phi$ in a slightly larger chart to obtain $\phi_m\to\phi$ uniformly, with Lipschitz constants at most $\ell+o(1)$ and upper Hessian bounds $D^2\phi_m\le(C_\phi+o(1))I$. The chain rule yields
\begin{align*}
 D_y^2[F(y,\phi_m(y))]
 &=F_{yy}+F_{yu}\otimes D\phi_m+D\phi_m\otimes F_{uy}
       +F_{uu}D\phi_m\otimes D\phi_m+F_uD^2\phi_m\\
 &\le\left(-\frac{c_0}T+2C(\ell+1)+CT(\ell+1)^2
                                 +e^{KT}(C_\phi+1)\right)I.
\end{align*}
For small $T$ the right side is at most $-c_0/(2T)I$. The semiconcavity inequality passes to the uniform limit. Positivity of $F_u$ is what permits use of the upper Hessian bound without any lower Hessian bound on the approximations. All maximizers lie in the smaller localization ball, so strict concavity on the larger convex ball proves global uniqueness.
\end{proof}

\begin{Pro}\label{prop:chain}
Let $f\in C^1(\R^m\times\R)$ satisfy $f_u>0$, and let $\phi$ be continuous. For $G(x)=f(x,\phi(x))$,
\[
 D^+G(x)=f_x(x,\phi(x))+f_u(x,\phi(x))D^+\phi(x).
\]
The formula is local and hence also applies in charts on $M$.
\end{Pro}
\begin{proof}
A supergradient can be represented by the derivative of a $C^1$ upper test function touching at the point. If $\eta$ is such a test for $\phi$, monotonicity makes $f(\cdot,\eta(\cdot))$ an upper test for $G$, proving one inclusion. Conversely, let $\zeta$ touch $G$ from above. The implicit function theorem solves $f(z,\eta(z))=\zeta(z)$ near $x$. Since $f_u>0$, this $\eta$ touches $\phi$ from above. Differentiating that identity gives the opposite inclusion. This argument needs only continuity of $\phi$. In the semiconcave applications all the superdifferentials are compact and nonempty.
\end{proof}

\section{Proofs of the calibration results}\label{app:B}
\subsection*{Proof of Proposition~\ref{prop:domination}}
Let $0<t_1<t_2$ and let $\xi$ be an admissible curve joining $x$ to $y$. Since $u$ is locally Lipschitz, the function $f(\tau)=u(\tau,\xi(\tau))$ is absolutely continuous. At almost every time, the semiconcave chain rule and the viscosity subsolution inequality give
\begin{equation}\label{eq:chain-domination}
 f'(\tau)\le L(\tau,\xi(\tau),\dot\xi(\tau),f(\tau)).
\end{equation}
Indeed, for each $(q,p)\in D^+u(\tau,\xi(\tau))$,
$q+H(\tau,\xi,p,f)\le0$, while Fenchel's inequality gives
$q+\langle p,\dot\xi\rangle\le L(\tau,\xi,\dot\xi,f)$.
The directional derivative formula \eqref{eq:SCfacts} gives \eqref{eq:chain-domination} at times where both compositions are differentiable.

Let $u_\xi$ solve \eqref{eq:Caratheodory} with initial value $f(t_1)$. Scalar comparison yields $f(t_2)\le u_\xi(t_2)$. Taking the infimum over $\xi$ proves \eqref{eq:domh}. Alternatively solve the scalar equation backward with terminal value $f(t_2)$. The same differential inequality, integrated backward, gives $u_\xi(t_1)\le f(t_1)$. Taking the supremum of these initial values proves \eqref{eq:dombh}.

At $t_1=0$ for $u=T_0^t\phi$, the first inequality is directly the variational formula. For the second, fix a curve from $x$ to $y$. Its endpoint map $E_\xi$ is strictly increasing and
$u(t_2,y)\le E_\xi(\phi(x))$. Thus
$E_\xi^{-1}(u(t_2,y))\le\phi(x)$. The supremum over curves gives \eqref{eq:dombh}. This avoids applying a chain rule to the possibly discontinuous initial datum.\qed

\subsection*{Proof of Proposition~\ref{prop:cal-equivalence}}
Set $f(s)=u(s,\eta(s))$ and let $U$ solve the scalar equation with $U(a)=f(a)$. The inequality \eqref{eq:chain-domination} gives $f\le U$. To see why endpoint equality propagates, define
\[
 k(s):=L(s,\eta(s),\dot\eta(s),f(s))-f'(s)\ge0.
\]
With $w=U-f$, the scalar Lipschitz bound gives
$w'=c(s)w+k(s)$, $|c(s)|\le K$, and hence
\[
 w(b)=\int_a^b\exp\!\left(\int_\sigma^bc(r)\,dr\right)k(\sigma)\,d\sigma.
\]
If (ii) holds, $w(b)=0$ forces $k=0$ almost everywhere. Uniqueness of the scalar equation yields $f=U$ on $[a,b]$. If (iii) holds, its difference from the integral of $f'$ is $\int_a^bk=0$, so the same conclusion follows. Conversely, $f=U$ gives both endpoint identities and their counterparts on every subinterval.

The curve is minimizing, since domination gives
\[
 u(b,\eta(b))\le u(a,\eta(a))+
            h(a,b,\eta(a),\eta(b),u(a,\eta(a)))\le U(b),
\]
and both endpoints of this chain are equal. Herglotz regularity therefore makes $\eta$ and its scalar lift $C^2$. It is calibrated in the sense of Definition~\ref{def:calibrated}, proving (i). The implication (i)$\Rightarrow$(ii) follows by taking the whole interval.\qed

\subsection*{Proof of Theorem~\ref{thm:cal-lifts}}
First consider an incoming calibrated curve $\gamma:[0,t]\to M$ and an interior time $0<s<t$. Choose nearby times $c<s<d$ on the curve. By the short-time estimates, after choosing $c,d$ sufficiently close to $s$, the endpoint functions
\begin{align*}
 (\sigma,z)&\longmapsto u(c,\gamma(c))+
             h(c,\sigma,\gamma(c),z,u(c,\gamma(c))),\\
 (\sigma,z)&\longmapsto u(d,\gamma(d))-
             \breve h(\sigma,d,z,\gamma(d),u(d,\gamma(d)))
\end{align*}
are smooth near $(s,\gamma(s))$. Domination puts $u$ below the first and above the second, and calibration gives equality there. The endpoint formulas \eqref{eq:negativevariation} and \eqref{eq:fullfirstvariation} give the same differential $(-H,p(s))$. Thus $u$ is differentiable at the interior point and its gradient has the asserted value. At the terminal time, choose \(c<t\) sufficiently close to \(t\). Then the function
\[
 z\longmapsto u(c,\gamma(c))
 +h(c,t,\gamma(c),z,u(c,\gamma(c)))
\]
is smooth near \(x=\gamma(t)\). By domination it lies above \(u(t,\cdot)\), while calibration gives equality at \(x\). Hence it is a smooth upper support for \(u(t,\cdot)\) at \(x\). By the endpoint formula \eqref{eq:negativevariation}, its gradient at \(x\) is \(p(t)\). Therefore $p(t)\in D_x^+u(t,x)$.

Now let $(q,p)\in D^*u(t,x)$. Choose differentiability points $(t_j,x_j)\to(t,x)$ with
$Du(t_j,x_j)\to(q,p)$, and take an incoming calibrated lift at each point. At a differentiability point its terminal momentum is $D_xu(t_j,x_j)$. Its time derivative is $-H(t_j,x_j,\allowbreak D_xu(t_j,x_j),\allowbreak u(t_j,x_j))$. Hence
$q=-H(t,x,p,u(t,x))$.
The a priori estimates and continuous dependence of the contact flow yield a limiting lift on $[0,t]$ with terminal data $(x,p,u(t,x))$.

There is a small initial-time issue because $\phi$ is only lower semicontinuous. If $\gamma_j$ are the approximating curves and $U$ is the limiting scalar lift, then
\[
 \phi(\gamma(0))\le\liminf_j\phi(\gamma_j(0))=U(0).
\]
Start a new scalar solution along the limiting curve at $\phi(\gamma(0))$. Its terminal value is at most $U(t)=u(t,x)$ by scalar comparison and at least $u(t,x)$ by the variational formula. Strict scalar order forces equality of the initial values. Consequently the limiting lift is calibrated, including at time zero, and is exactly \eqref{eq:reachableflow}.

Conversely, approach the terminal point of any calibrated lift by its interior points. The differentiability just proved and continuity of the lift give
\[
 Du(s,\gamma(s))\longrightarrow
 (-H(t,x,p(t),u(t,x)),p(t))\qquad(s\uparrow t).
\]
This is a reachable space-time gradient. Finally, uniqueness of the contact ODE identifies any two lifts having the same terminal state, proving injectivity and the desired bijection.\qed

\bibliographystyle{plain}
\bibliography{mybib}
\end{document}